\magnification=\magstep1
\input amstex
\documentstyle{amsppt}

\def\bbZ{{\Bbb Z}}

\def\bbC{{\Bbb C}}
\def\bbQ{{\Bbb Q}}
\def\bbP{{\Bbb P}}

\def\bbF{{\Bbb F}}
\def\bfF{{\bold F}}

\def\calO{{\Cal O}}

\def\calE{{\Cal E}}

\def\Pic{\roman{Pic}}

\hsize 6in

\topmatter

\title 
COBLE RATIONAL SURFACES
\endtitle

\author 
Igor V. Dolgachev and De-Qi Zhang
\endauthor

\thanks 
Research of the first author was supported in part by a NSF grant 
\endthanks

\endtopmatter
\document

\head Introduction \endhead

A Coble surface is a nonsingular projective rational surface $S$ with empty 
anticanonical linear system $|-K_S|$ but non-empty bi-anticanonical system
$|-2K_S|$. A classical example of such surface is the blow-up of $\bbP^2$ at 10
nodes of  an irreducible rational plane curve of degree 6 with ordinary  nodes as
singularities. Rational plane sextics of this kind were intensively  studied by A.
Coble [Co1,Co2]. Among other things he showed that a Cremona  equivalence classe of
such curve is the union of finitely many projective equivalence classes (see
[Co2],[MS]). This result can be interpreted as saying that the automorphism group of
the associated Coble surface is isomorphic to a subgroup of finite index in the
orthogonal group of the lattice
$M_S = (K_S)^\perp_{\roman{Pic}(S)}$. In fact, Coble shows that for a general
sextic with 10 nodes this group is isomorphic
to the congruence subgroup of level 2 of the
group
$O(M_S)/\{\pm 1\}$ (see [Do2]). The latter group is isomorphic to the Weyl
group of infinite root system of type $E_{10}$.  A similar answer is
known for a generic Enriques surface. It was obtained much later by V.
Nikulin [Ni] and independently W. Barth and C. Peters [BP]. 

The connection between classical Coble surfaces and Enriques surfaces is a
nice one: a double cover of a Coble surface branched along the proper
transform of a reduced sextic is a K3 surface which is a degeneration of  the
K3-cover of an Enriques surface. More geometrically, the embedding of a
Coble surface in  $\bbP^5$ defined by the linear system of curves 
of degree 10 with singular points of multiplicity 3 at the 10 nodes 
of the sextic is a surface of degree 10 which is a projective degeneration 
of an Enriques surface in its Fano embedding in  $\bbP^5$. 

Another result of Coble is that the set of rational 
smooth curves with negative self-intersection on the blow-up of
ten nodes of an irreducible sextic is finite modulo the automorphism group
[Co2]. It is known that the same fact is true for all minimal
non-rational  algebraic surfaces [DP] and it is not true for any general
blow-up of  $\ge 9$ points in
$\bbP^2$. This makes some (not all, as we shall see) of Coble surfaces
exceptional in this respect. In fact,  this work in which we classify all
Coble surfaces was partially motivated by the problem to classify all
rational surfaces with finitely many smooth rational curves of negative
self-intersection modulo automorphisms of the surface.

The classification of Coble surfaces is related to other classification problems. 
To be precise, when $|-2K_S|$ contains a reduced divisor, the double cover branched
along this divisor is a normal K3 surface with an involution. The
classification of K3-surfaces with an involution can be found in [Zh3]
extending some earlier results of Nikulin [Ni]. 
In particular, a terminal Coble surface is the minimal resolution of
a maximum rational log Enriques surface of
index 2 in the sense of [Zh2, OZ] (cf. Proposition 6.4).  
The latter surfaces were classified in [Zh1, Zh2]. 

Let us describe the main results of this paper.  First of all we divide Coble surfaces 
into two major classes. For a surface of the first class ({\it elliptic type}; cf. 2.9) 
there exists a birational morphism onto a surface $Y$ such that the anticanonical 
linear system $|-K_Y|$ has only one member, and a general member of the mobile part of
$|-2K_Y|$ is a smooth elliptic curve. Surfaces of the second class ({\it rational
type}) admit a similar birational morphism only this time the mobile part
of  $|-2K_Y|$ consists of divisors of arithmetic genus 0 (not necessarily
irreducible). We show that Coble surfaces of elliptic type are
obtained as either blow-ups of singular points 
and their infinitely near points of a non-multiple fibre on
a minimal rational elliptic surface with one multiple fibre of
multiplicity 2 ({\it Halphen type}), or as blow-downs of some disjoint
sections and maybe components of one fibre of a non-minimal rational
elliptic surface with a section ({\it Jacobian type}). We also give a
construction for surfaces of rational type as blow-ups of minimal
rational surfaces. It turns out that surfaces of elliptic type always
admit a birational morphism to $\bbP^2$.  However, for any given  $n$
there are Coble surfaces of rational type which do not admit a birational 
morphism to the minimal ruled surface  $\bold F_n$. We prove that Coble surfaces 
of elliptic type are obtained by blowing up  $\bbP^2$  with centers at 
singular points of certain plane curves  $\Gamma$  of degree 6 
(Coble sextics), but the center of the very last blow up may not be on 
$\Gamma$.  We describe such sextics. 

An important class of Coble surfaces $X$ to which the original example of 
Coble belongs is the one where the linear system $|-2K_X|$ contains 
a reduced divisor.  In this case $X$ admits a double cover 
which is a K3-surface with at most ordinary double points. 
We prove, under an appropriate condition of generality, 
that surfaces of this kind contain only finitely many smooth rational curves 
with negative intersection modulo automorphisms of the surface. 
We show in Example 6.10, that this statement cannot be extended to all Coble surfaces. 
It is possible that every Coble surface of rational type contains only 
finitely many negative rational curves; however we could not prove it.

Finally a word about the ground field $k$. We assume it to be algebraically
closed, though most of the paper does not use any assumption on the characteristic.
However Theorem 6.7 assumes that $k = \bbC$ (with more efforts one can give 
another proof which does not use this assumption) and we assume that $k$ is
uncountable in Example 6.10.

\medskip
{\bf Acknowledgement.} This joint work was done during the second 
author's visit at University of Michigan in Summer of 1999, who would
like to thank its Department of Mathematics for the hospitality. The
first author would like to acknowledge the support and hospitality of the
Mathematical Sciences Research Institute at Berkeley during June of 1999.

\medskip
\head 1. Some preliminary results \endhead

\noindent
{\bf 1.1} We shall consider an order on the set of Coble surfaces defined by dominant
birational morphisms $f:X'\to X$. Thus we can speak about a {\it minimal Coble 
surface} $X$ (which does not admit a birational, but not biregular, morphism onto another 
Coble surface) and a
{\it terminal Coble surface} which is not the image of any birational 
but not biregular morphism of Coble surfaces. 

We shall see that 
there exist minimal and terminal Coble surfaces, as well
as non-minimal or non-terminal Coble surfaces (Example 2.6).
It follows from the Riemann-Roch theorem that $K_X^2\le -1$ for any Coble surface. 
Hence a Coble surface with $K_X^2 = -1$ is always minimal. We shall see 
that there are also minimal Coble surfaces with $K_X^2 < -1$ (Example 4.9).

In the next paragraphs we shall give some conditions for a Coble surface to be minimal.

\medskip\noindent
{\bf 1.2} 
For any positive divisor $D$ on a nonsingular projective surface $V$
we set
$$p_a(D) = {1\over 2}(D^2+K_V\cdot D) +h^0(\calO_D) = h^1(\calO_D).\eqno (1.1)$$
Here the last equality follows from Riemann-Roch, applied to the divisor $D+K_V$.

\medskip
\plainproclaim Lemma. Assume $h^1(\calO_V) = h^2(\calO_V) = 0$, for example,
$V$  is a rational surface. Then we have:
\roster
\item
$p_a(D) = h^0(D+K_V)$; in particular, $p_a(D_1) \le p_a(D)$  if  $D_1 \le D$.
\item
If  $p_a(D_1) \ge 1$ (e.g., when the dual graph of  $D_1$  contains a loop),
then  $p_a(D_1+D_2) \ge h^0(D_2)$.
\endroster

{\sl Proof.} We only need to show the first part of (1).  
It follows from considering the exact sequence
$$0\to \calO_V(-D)\to \calO_V \to \calO_D \to 0.$$
Indeed, the sequence implies that 
$H^1(D,\calO_D) \cong H^2(V,\calO_V(-D)) = H^0(V,\calO(D+K_V)$. 

\medskip\noindent
{\bf 1.3}
The part $h^0(\calO_D)$ in (1.1) is often hard to compute. We only cite the
following useful results which can be found for example in [Re], p.81. Recall
that an effective divisor is called {\it numerically $k$-connected}, 
if for any decomposition $D = D_1+D_2$
into positive parts, we have $D_1\cdot D_2 \ge k$.

\medskip
\plainproclaim Lemma.
\roster 
\item Let $0 < D' < D$. Then there is an exact sequence 
$$0\to \calO_{D-D'}(-D')\to \calO_D\to \calO_{D'}\to 0.$$
\item
Assume that  $D$  is  numerically 1-connected. Then
$$h^0(\calO_D) = 1.$$
\item In particular, let  $D_i > 0$  such that  $D_1 + D_2$ 
is reduced and  $D_1, D_2, D_1+D_2$  are all numerically 1-connected.
Then  
$$p_a(D_1+D_2) = p_a(D_1)+p_a(D_2) + D_1\cdot D_2 - 1.$$
\endroster

\noindent \medskip   
\plainproclaim Lemma 1.4. Let  $X$  be a Coble surface.  Then any member
$D$  of  $|-2K_X|$  consists of smooth rational curves and is of simple normal crossing.

{\sl Proof.} If Lemma 1.4 is false, then
$D$  contains a reduced connected divisor  $D_1$  such that either
$D_1$  is irreducible with arithmetic genus $\ge 1$,
or  $D_1$  is a loop, or  $D_1$  is the sum of two 
curves with an order $\ge 2$  contact at a point,
or  $D_1$  is the sum of 3 curves sharing one point.
This leads to  $1 \le p_a(D_1) \le p_a(D) = h^0(-K_X) = 0$ (Lemma 1.2),
a contradiction.  Hence Lemma 1.4 is true.

\medskip
By an {\it exceptional curve} we shall mean a one-dimensional fibre of a birational
morphism between nonsingular projective surfaces. An irreducible exceptional curve
is a $(-1)$-curve. Here by a {\it $(-n)$-curve} we mean a smooth rational curve $R$ with 
$R^2 = -n < 0$.

\medskip
\plainproclaim 1.5 Lemma. Let  $\pi : X \rightarrow Y$
be the blow-down of a  $(-1)$-curve  $E$  on a smooth rational
surface  $X$  to a point  $q$  on  $Y$.
\roster
\item
Suppose that  $X$  is a Coble surface.
Then for any $D\in |-2K_X|$ and any  $s \ge 1$, one has
$$p_a(D+sE) = h^0(-K_Y).$$
In particular, $p_a(D+2E)\le 1$; $Y$ is also a Coble surface if and only if 
$p_a(D+2E) =0$.
\item
Suppose that  $h^0(-K_Y) = 1$.
Then  $X$  is a Coble surface if and only if  $q$  is a multiplicity
$\ge 2$  point of a member in  $|-2K_Y|$  but  $q$  is not a point
of the unique member in  $|-K_Y|$.
\item
Suppose that  $Y$  is a Coble surface.  Then  $X$  is also a Coble
surface if and only if  $\text{\rm mult}_q(D) \ge 2$  for some  $D \in |-2K_Y|$.  
\endroster

{\sl Proof.}  It follows from the projection formula and Lemma 1.2 that
$$h^0(Y,-K_Y) = h^0(X, \pi^*(-K_Y) + (s-1)E) = h^0(X,-K_X + sE) = h^0(X, K_X + D + sE)
= p_a(D + sE).$$

\par
To prove the last part in (1),
we first note that $|-2K_Y| = |\pi_*(-2K_X)|\ne \emptyset$, so that $Y$ is a Coble
surface if and only if $p_a(D+2E) = 0$. Also, notice that $h^0(-K_Y) \le 1$
since otherwise we can find an anticanonical divisor on $Y$ which passes through
the point $\pi(E)$. This would imply that $|-K_X|\ne \emptyset$ contradicting the
assumption that $X$ is a Coble surface.

For (2) and (3), see Lemma 1.9 below.

\medskip
\plainproclaim 1.6 Corollary. A Coble surface $X$ is minimal if and
only  if, for any
$(-1)$-curve $E$ and any $D\in |-2K_X|$, 
$$p_a(D+2E) = 1.$$

\medskip
The next lemmas will be used frequently in the subsequent sections.

\medskip
\plainproclaim  1.7 Lemma. Let  $X$  be a smooth rational surface.
Then we have:
\roster
\item  Suppose that  $L$  is a smooth rational curve with  $L^2 \ge 0$.
Then  $|L|$  is base point free and  $h^0(Y,L) = 2 + L^2$.
\item  Suppose that $L$  is an irreducible curve with  $p_a(L) = 1$
and  $L^2 \ge 1$.  Then a general member of  $|L|$  is
smooth and  $h^0(Y,L) = L^2 + 1$.  
If  $L^2 \ge 2$, then  $Bs|L| = \emptyset$.
If  $L^2 = 1$  then  $|L|$  has exactly one base point.
\item Suppose that  $L$  is smooth elliptic with  $L^2 = 1$
and  $G$  an effective divisor, not linearly equivalent to  $L$,
such that  $-2K_X \sim L + G$ (if  $h^0(X, -K_X) \le 1$  
one always has  $G \notin |L|$).  
Let  $G_1$  be the unique component in  $G$  with  $L\cdot G_1 = 1$.
Then  $G_1 \notin |L|$  and
$L \cap G_1 = L \cap G$  is the unique base point of  $|L|$.
\endroster

{\sl Proof}. (1) follows from the exact sequence below,
the induction on  $L^2$  and the fact that the result is true
when  $L^2 = 0$:
$$0 \rightarrow \calO_X \rightarrow \calO_X(L) \rightarrow
\calO_{\bbP^1}(L^2) \rightarrow 0.$$

For (2), a similar exact sequence as in (1) shows that
$h^0(X, L) = 1 + h^0(L, \calO_L(L)) = 
1 + h^0(L,\omega_L(-L)) + 1 - p_a(L) + L^2 = 1 + L^2$.
Here  $\omega_L \cong {\Cal O}_L$  is the dualizing sheaf and
we have applied the duality and the Riemann Roch theorem 
for  $L$.

We assert that a general member of  $|L|$  is
smooth.  Take $L^2-1$ generic points such that the linear system of divisors from
$|L|$ passing through these points is one-dimensional. By blowing up the points on 
$L$, it suffices to prove the assertion when  $L^2 = 1$. Since the only base
point of $|L|$ is then simple, the assertion follows from Bertini's theorem. 

For (3), if  $G_0 \le G$  and  $G_0 \in |L|$, then 
$0 \sim (K_X + L) + (K_X + G_0) + (G - G_0) \ge G - G_0$,
which leads to  $G = G_0 \sim -K_X \sim L$, a contradiction.
It remains to show that  $L \cap G_1$  is equal to
the unique base point  $p$  of  $|L|$.
Let  $\tau : Y \rightarrow X$  be the blow-up of
$p$  with  $C$  the  $\tau$-exceptional curve.
Then  $-2K_Y \sim \tau^{-1}(L) + \tau^{-1}(G) - C$ (we use $\tau^{-1}$ 
to denote the proper inverse transform
under a birational map).  Since $-2K_Y$, $\tau^{-1}(L),\tau^{-1}(G-G_1)$ 
can be represented by a divisor 
contained in fibres, we obtain that the restriction of $\tau^{-1}(G_1)$ 
and $C$ to a general fibre is linearly
equivalent. This is obviously impossible (since no two distinct points on 
an irrational curve are linearly equivalent). 

\noindent
\plainproclaim 1.8 Lemma. (M. Miyanishi) Let  $V \to \bbP^1$  be a smooth 
rational ruled surface with two sections  $s_1, s_2$.  Then there is a birational morphism
$\pi: V \to {\bold F}_d$  onto a minimal ruled surface of degree  $d$, such that
$\pi(s_1)\cdot \pi(s_2) = s_1\cdot s_2$.  Moreover, if both  $s_i^2$  are 
negative we can choose  $\pi$  such that  $\pi(s_1)^2 = -d = -1$.
Finally, $\pi(s_2)^2 = -\pi(s_1)^2 + 2(s_1 \cdot s_2)$.

{\sl Proof.} Let 
$\pi_1 : V \rightarrow V_1$  be the composition of smooth blow-downs of all  
$(-1)$-curves in fibers disjoint from  $s_1$
and $s_2$. Since both  $s_i$  are sections, we see easily that
for every singular fiber  $F_i$  on  $V_1$, $s_1 + F_i + s_2$  has the following 
dual graph:
$$s_1 - (-1) - (-2) - \cdots - (-2) - (-1) - s_2.$$
Note that  $\pi_1(s_i)\cdot \pi_1(s_j) = s_i\cdot s_j$.
Now a suitable blow-downs of  $(-1)$-curves in fibers on
$V_1$  will give the required birational morphism  $\pi$.
For the second assertion, we let  $V_1 \rightarrow {\bold F}_1$
be the successive blow-downs of  $(-1)$-curves in fibers
with exactly  $-1 - s_1^2$  of them intersecting  $s_1$;
this is possible because a minimal ruled surface has at most
one negative section.  The last assertion follows by expressing  $\pi(s_2)
\sim \pi(s_1) + [s_1 \cdot s_2 - \pi(s_1)^2] f$, where $f$ denotes a fibre.

\noindent
\plainproclaim 1.9 Lemma. Let  $X_1 \rightarrow \cdots \rightarrow X_n$
($n \ge 2$) be a sequence of blow-ups  
$\tau_i : X_i \rightarrow X_{i+1}$  of smooth surfaces 
with center  $p_{i+1} \in X_{i+1}$  and exceptional curve  $E_i \subset X_i$.
\roster
\item 
Assume that a positive divisor  $D$  belongs to  $|-2K_{X_1}|$
and denote by  $D_i$  its  direct image on  $X_{i+1}$.
Then each  $p_{i+1}$  is a singularity of  $D_{i+1}$.
In particular, the divisor  $D_n$  is always singular.
\item
If  $|-2K_{X_1}| \ne \emptyset$, then there is a singular member  $D_n$
of  $|-2K_{X_n}|$  such that the indeterminancy locus
of the rational map  $X_n \cdots \rightarrow X_1$  is contained in
the singular locus of  $D_n$.
\endroster

{\sl Proof.} This follows from the fact that  $D_i \in |-2K_{X_i}|$
and hence  $E_i\cdot D_i = 2$.

\medskip
In view of the next result, which follows from the fact that
a Coble surface  $X$  always has  $K_X^2 \le -1$,
we only have to consider minimal Coble surfaces.

\medskip \noindent
\plainproclaim 1.10 Lemma. 
Suppose that  $X$  is a Coble surface.  Then there is a sequence of
blow-downs  $X = X_1 \rightarrow X_2 \rightarrow \cdots \rightarrow
X_n$ ($n \ge 2$)  such that  $X_n$  is not Coble but  $X_i$ ($i < n$) 
are all Coble and especially, $X_{n-1}$  is a minimal Coble
(see 1.5 and 1.9 for the restriction on the centers of blow-ups).

\head 2. The elliptic case \endhead

\medskip \noindent
{\bf 2.1} Let  $X$  be a Coble surface and $E$ a $(-1)$-curve
with  $\pi : X \rightarrow Y$  the blow-down of  $E$.
Assume the hypothesis that  $p_a(-2K_X+2E) = 1$, i.e., $|-K_Y| \ne \emptyset$
(on a minimal Coble surface, any  $(-1)$-curve satisfies this, by Lemma 1.5).
  
Consider the linear system
$$|-2K_X+2E| = \pi^*(|-2K_Y|).$$
Note that  $\dim |-2K_Y|> 0$  for otherwise  $|-2K_Y| = 2|-K_Y|$  and 
hence  $|-K_X| \ne \emptyset$, a contradiction (cf. Lemma 2.3 below).

Write
$$|-2K_X+2E| = |M| +P, \eqno (2.1)$$
where $|M|$ is the mobile part, and $P$ the fixed part.
We also write  
$$P = G + H, \,\,\,\, G = \sum_{i=1}^J g_i G_i,$$ 
where  $G_i\cdot M \ge 1$  while  $H\cdot M = 0$.

\par
By Lemmas 1.2 and 1.5, $p_a(M) \le p_a(M+P) = 1$.
We say that $X$ is of {\it elliptic type} with respect to  $E$
if  $p_a(M) = 1$, and of {\it rational type} with respect to  $E$ if 
$p_a(M) = 0$.  It may happen that the same Coble surface
$X$ (even minimal one) is of elliptic type with respect to one  $E_1$, but of
rational type with respect to another  $E_2$ (see Example 2.11).  
The  $(-1)$-curves here are used like markings to help classify Coble surfaces. 
A minimal Coble surface will be called of rational type 
if it is of rational type with respect to any  $E$.

\medskip
\plainproclaim 2.2 Lemma. Let  $X$  be a Coble surface with
a  $(-1)$-curve  $E$  satisfying  $p_a(-2K_X+2E) = 1$
(on a minimal Coble surface, this is always true 
for any  $(-1)$-curve, by Lemma 1.5).
Then, with the above notation, we have the following:
\roster
\item
If  $p_a(M) = 1$, then a general member of
$|M|$  is a smooth elliptic curve, and  $G \cdot M = M^2$.
\item
If  $p_a(M) = 0$, then  either  $M = kM_1$ ($k \ge 1$)
with  $M_1 \cong \bbP^1$, $M_1^2 = 0$  and  $G \cdot M_1 = 4$, or
$M \cong \bbP^1$, $M^2 \ge 1$  and  $G \cdot M = 4 + M^2$.
\newline
In the following, we set  $M_1 = M$, when  $M$  is
irreducible.
\item
$M_1 + P$  is of simple normal crossing.
$P$  consisits of  $(-n)$-curves  with  $n \ge 1$.
\item 
Suppose that  $D_1 + \cdots + D_s$  is a chain
in  $P_{\text{\rm red}}$, such that
$L := M_1 + \sum_i D_i$  is a loop.
Then  $p_a(M) = 0$, $L$  is a simple loop, and
$\sum_i D_i^2 \le -2s - 1$; moreover, $L$  is the only
loop in  $M_1 + P_{\text{\rm red}}$.
\endroster

{\sl Proof.} For (1) and (2), we have only to consider the case
where a general member of  $|M|$  is reducible.
Then by Stein factorization and the rationality of  $X$
(or rather the vanishing of  $q(X)$), we have
$M = kM_1$ ($k \ge 2$) with  $|M_1|$  an irreducible pencil.
Since  $p_a(M_1) \le p_a(M) \le 1$, Lemma 1.7 and the fact that
$\dim |M_1| = 1$ imply that either  $p_a(M_1) = 0$  and
$M_1^2 = 0$, or  $p_a(M_1) = 1$  and  $M_1^2 \le 1$.
The latter case leads to that  $1 \ge p_a(M) \ge p_a(2M_1) \ge 2$
(see Lemma 1.3 when  $M_1^2 = 1$), a contradiction.
This proves (1) and (2); indeed the equality on
$G \cdot M$  or  $G \cdot M_1$  is obtained by intersecting
both sides of (2.1) with  $M_1$.

Next we prove (3).  First  $P$  is of simple normal crossing
and contains no arithmetic genus  $\ge 1$  curves,
for otherwise, $1 \ge p_a(M_1+P) \ge h^0(M_1) \ge 2$
by Lemma 1.2.  Thus each curve in  $P$  is a
$(-n)$-curve with  $n \ge 1$  because  $P$
is the fixed part and by Lemma 1.7.
If  Bs$|M| = \emptyset$  then (3) is clear.
So, in view of Lemma 1.7, we only need to consider
the case where  $p_a(M) = 1$.  Then (3) can be proved
in a manner similar to (4) below.

Now we prove (4).  If  $p_a(M_1) = 1$  or the dual graph of  $L$  is not
a simple loop then  $2 \le p_a(L) \le p_a(M+P) \le 1$;
if the dual graph of  $M_1+P_{\text{\rm red}}$  contains 
another loop, then there is a linear chain  $N$  having no common components 
with  $L$  such that  $N \cdot L \ge 2$, which leads to
$2 \le p_a(L+N) \le p_a(M+P) \le 1$,
again a contradiction.  As in Lemma 1.7, one sees easily that  
$h^0(L) = 1 + h^0(L|L) \ge 1 + \chi(L|L) =
1 + \chi({\Cal O}_L) + L^2 = 1 + L^2$.
Substituting  $h^0(L) = h^0(M_1) = M_1^2 + 2$  and expanding
$L^2$, we will get the inequaltiy in (4).

\medskip
\plainproclaim 2.3 Lemma. Let $X, E$  and notation be as in Lemma 2.2.
Then we have
\roster
\item
$E \cap P = \emptyset$, whence  $E\cdot M = E \cdot P = 0$,
\item
$|M|$  and  $P$  are pull backs of
$|\pi(M)|$  and  $\pi(P)$, whence a general member
$M$  is disjoint from  $E$, and
\item
$|M|$ contains a member  $M'$  with  $M'-2E \ge 0$.
\endroster

{\sl Proof.} Suppose $E\cap P\ne \emptyset$. Then the linear system 
$|-2K_Y|$  has the point $p = \pi(E)$ as a base point. 
Let $C$ be the unique divisor in $|-K_Y|$. The
divisor $2C\in |-2K_Y|$ and hence contains $p$. Thus $p\in C$ and hence $|-K_X|\ne
\emptyset$. This contradiction proves the first assertion, which, in turn,
implies the rest.

\medskip
In the rest of the section, we shall classify
Coble surfaces  $X$  of elliptic type (see 2.5, 2.6, 2.8, 2.9).

\medskip\noindent
{\bf 2.4} Case: $p_a(M) = 1$  and  $M^2 = 0$. 
In this case, $|M|$ is a pencil of elliptic curves without base-points.
It defines an elliptic fibration $f:X\to \bbP^1$, so that  
$E$  and  $P$  are contained in fibres (Lemma 2.3).
Blowing down $E$, we get an elliptic fibration
$f_Y: Y\to \bbP^1$. Let $f_m : Y_m \to \bbP^1$ be its relative minimal model, i.e.
$Y_m$ is obtained from $Y$ by blowing down exceptional curves contained in fibres of
$f_Y$. 

Recall that a relative minimal rational elliptic surface $V$ is called an {\it Halphen
surface of index $n$} if the divisor class of its fibre is equal to $-nK_V$. Any
relative minimal rational elliptic surface is an Halphen surface of some index $n$. 
An Halphen surface of index 1 is a {\it Jacobian} rational elliptic surface. 
It is characterized by the condition that 
the fibration does not have multiple fibres, or
equivalently, admits a section. An Halphen surface of index $n \ge 2$ has a unique
multiple fibre $nF_1$ of multiplicity $n$. In this case 
$$|-iK_V| = \{iF_1\}, \quad
1\le i \le n-1 .\eqno (2.2)$$
All of this is rather well-known and can be found for example in [CD], Chapter 5, \S 6.

Let $n$ be the index of $f_m : Y_m \to \bbP^1$. Since $|-K_Y|\ne
\emptyset$ and by Lemma 1.9, $Y$ is obtained from $Y_m$ by a successive blow-ups 
of {\it singular} points and their infinitely near points on {\it one} fibre $F_1$ 
(the unique multiple fibre if $n \ge 2$). We know
that $\dim |-2K_{Y_m}| \ge \dim |-2K_Y| \ge 1$. 
Applying (2.2) this implies that $n \le 2$. Moreover, if $n =2$, after one
blow-up the anti-bicanonical linear system becomes of dimension 0. So, in this case,
we must have $Y = Y_m$ and $P = 0$. 

Suppose that  $n = 1$.  We claim that  $X$  is not a minimal Coble surface.
Since  $h^0(Y_m, -2K_{Y_m}) = 2$  while  $h^0(Y, -K_Y) = 1$ (Lemma 1.5),
we have  $Y \ne Y_m$.  For simplicity, we assume that
$Y \rightarrow Y_m$  is a single blow-down of a  $(-1)$-curve  $E_1$  to
a point  $q_1$  on a fiber  $F_1$  (the general case is similar).
Then the mobile part of  $|-2K_Y|$  is equal to the pull back of 
the elliptic pencil, and its fixed part is equal to [(the proper inverse transform
of  $F_1$) $+ (m_1-2)E_1$], where
$m_1$  is the multiplicity of  $F_1$  at  $q_1$.
The map  $\pi : X \rightarrow Y$  in  2.1
is just the blow-down of the  $(-1)$-curve 
$E$  to a multiplicity  $m$ ($\ge 2$) singular point of a fibre  
$F \ne F_1$ (cf. Lemma 2.3 and Remark 2.9 below).  Moreover,
$$|M| + P = \pi^*|-2K_Y| = |F| + (F_1' + (m_1 - 2)E_1'),$$
$$-2K_X \sim P_1 + P_2, \,\,\, P_1 := F' + (m-2)E, \,\,\,
P_2 := F_1' + (m_1 - 2)E_1'.$$
Here  $F_1', F', E_1'$  denote the proper inverse transforms
of  $F_1, F, E_1$, and  $F$  denotes a full fiber on  $X$
by abuse of notation.

\par
Let  $C_m$  be a section on  $Y_m$  and  $C$
be its total inverse transform on  $X$.
Since we blow-up singular points of $F_1, F$, the
map  $X \to Y_m$ is an isomorphism over $C_m$.  Therefore $C$ is a $(-1)$-curve.
Let us check  that $p_a(-2K_X+2C) = 0$ so that after blowing down $C$ 
we get a Coble surface again.  Applying the exact sequence in Lemma 1.3 to
compare  $p_a(P_1+P_2+C)$  with  $p_a(P_1+P_2+2C)$ and 
$p_a(P_1+P_2)$, we find that  $p_a(P_1+P_2+2C) = p_a(P_1+P_2+C) = p_a(P_1+P_2)$. 
The latter equality follows from the fact that $h^0(\calO_{P_1+P_2}) > 
h^0(\calO_{P_1+P_2+S})$  and the application of Lemma 1.3 with  
$D - D' = C$.
Since $p_a(P_1+P_2) = h^0(K_X+P_1+P_2) = 0$, the claim is proved.

\par
Summing up, we obtain
\noindent
\plainproclaim 2.5 Theorem. 
Let  $X$  be a minimal Coble surface and  $E$ a $(-1)$-curve on $X$.
Assume that the mobile part  $|M|$  of $|-2K_X+2E|$ satisfies 
$p_a(M) = 1$  and  $M^2 = 0$.
Then $|-2K_X+2E| = |M|$ and  $X$ is obtained from an 
Halphen surface  $Y$  of index 2 by one blow-up  $\pi$
of a singular point on a non-multiple fibre $F$
with  $E$  the exceptional curve.  

\medskip\noindent
{\bf 2.6} {\bf Definition, Remark and Example}. 
(1) A Coble surface  $X$  is of {\it Halphen type}, or type(H), with respect  to 
$E$  if it is obtained as in Theorem 2.5 above.
In general, a Coble surface  $W$  is of {\it Halphen type}
if there is a birational morphism  $W \rightarrow X$
such that  $X$  is of Halphen type with respect to some  $E$.
 
(2) From 2.4 and 2.7 below, we see that an arbitrary Coble surface  $X$
with a  $(-1)$-curve  $E$  satisfying  $p_a(-2K_X+2E) = 1$,
$p_a(M) = 1$  and  $M^2 = 0$, is equal to either  $X$
in Theorem 2.5, or  $X'$  in Theorem 2.8 where  $q$  is
a singular point of  $F$ (cf. Remark 2.9).

(3) Let $Y\to \bbP^1$ be an Halphen surface and $F$ a
non-multiple singular fibre. If $F$ is of type $I_n$ (= $\tilde A_{n-1}$ 
in other notation), then $F$ has exactly  $n$ double points.
Blowing up one double point gives us a minimal Coble surface  $X$  because  $K_X^2 = -1$.
Blowing up all double points gives us a terminal 
(but non-minimal if  $n \ge 2$) Coble surface (cf. Proposition 6.4 in \S 6).
In particular, if $n = 1$, we get a Coble surface which is both minimal and terminal.
The same is true when we blow up the unique singular point of a fibre of type
$II$.  On the other hand, if we blow-up successively points on multiple components of
a non-reduced fibre, we get examples of non-terminal Coble surfaces.

\medskip\noindent
{\bf 2.7} Case :  $p_a(M) = 1$  and  $M^2 = m \ge 1$.
In notation of Lemma 2.2, we have  $M \cdot G_i = 1$  and
$G_i \cap G_j = \emptyset$  when  $i \ne j$.  
Thus  $\sum_{i=1}^J g_i = m$  by Lemma 2.2.

\par
By Lemma 1.7, $\dim |M| = m$. 
Fix a general member  $M_1$  of  $|M|$  and put  $p_i = M_1\cap G_i$.
Blow up the points $p_i$ to get a surface  $X_1$. 
If  $g_i \ge 2$, we pick the point  $p_i^{(1)}$  on
$X_1$ lying over  $p_i$ and on the proper inverse transform of $M_1$. 
Continue in this way to get a surface  $X'$  obtained from  $X$ 
by blowing up the points $p_i = p_i^{(0)}, p_i^{(1)},\ldots, p_i^{(g_i-1)}$, 
where  $p_i^{(s)}$  is infinitely near to  $p_i^{(s-1)}$  and
lies on the proper inverse transform of  $M_1$.  

\par
Let  $\Theta_i^{(j)}$  be the proper inverse transform on  $X'$
of the  $(-1)$-curve lying over the point  $p_i^{(j-1)}$.
Set  $C_i := \Theta_i^{(g_i)}$  and denote by
$M_1'$, $G_i'$ the proper inverse transforms on  $X'$  of  $M_1$, $G_i$.  
Then  $M_1' + C_i + \Theta_i^{(g_i-1)} + \cdots + \Theta_i^{(1)} + G_i'$
has the dual graph:
$$(0) - (-1) - (-2) - \cdots - (-2) - G_i'.  \eqno (2.3)$$

\par
Let $E'$ be the pre-image
of $E$ on $X'$.  Since each $p_i$ is not on $E$ (Lemma 2.3),
$X' \rightarrow X$  is an isomorphism over  $E$  and hence
$E'$ is a $(-1)$-curve.  Let  $\pi' : X' \rightarrow Y'$  
be the blow-down of  $E'$  to a point  $q$.
Then there is a birational morphism  $Y' \rightarrow Y$
such that two compositions  
$X' \overset{\pi'}\to{\rightarrow} Y' \rightarrow Y$
and  $X' \rightarrow X \overset{\pi}\to{\rightarrow} Y$  are identical.
Applying Lemma 1.7 to  $S$, which is the blow-down of  $C_1$  
to the point  $p_1^{(g_1-1)}$  followed by the blow-down of  $E'$,
and  $L :=$ (the image on  $S$  of  $M_1'$),
we obtain  $h^0(X', M_1') = h^0(S, L) = 2$.

\par
Noting that each
$p_i$ is a point of multiplicity  $1+g_i \ge 2$ in 
$M_1 + \sum_{i=1}^J g_iG_i + H$ ($\sim -2K_X + 2E$), 
we get
$$(\pi')^*(-2K_{Y'}) = -2K_{X'}+2E' \sim M_1' + P', \eqno (2.4)$$
where  $P'$  is the sum of the total transform of  
$H$  and the disjoint union of  $J$  weighted linear chains
$\Theta_i + g_iG_i'$  with  
$\Theta_i = \sum_{j=1}^{(g_i-1)}(g_i-j)\Theta_i^{(j)}$.

\par
$|M_1'|$  defines a Jacobian elliptic fibration with
sections  $C_i$.  Since  $P'\cdot M_1' = 0$  and  
$E'\cdot M_1' = E\cdot M = 0$
(Lemma 2.3), $P', E'$  are contained in fibres
$F_1, F$  on  $X'$, respectively.
Let  $Y' \rightarrow Y_{\min}$  be the smooth blow-down
to a relative minimal model.  As we explained in 2.4, 
$X'$  is obtained from  $Y_{\min}$  by blowing up singular points  
and their infinitely near points
on a fibre  $F_1$  of the elliptic fibration on  $Y_{\min}$  followed 
by blowing up a point  $q$  of another fibre  $F$  to the curve  $E'$.
Let us sum up the previous arguments by stating the following:

\medskip
\plainproclaim 2.8 Theorem. Let $X$ be a Coble surface with a $(-1)$-curve $E$
satisfying  $p_a(-2K_X+2E) = 1$ (on a minimal Coble surface,
any  $(-1)$-curve satisfies this).
Assume that the mobile part $|M|$ of $|-2K_X+2E|$ 
satisfies  $p_a(M) = 1$  and  $M^2 = m > 0$. Then $X$ is obtained as follows. 
\newline
There exist a relative minimal Jacobian rational elliptic surface 
$Y_{\min}$  with a singular fibre  $F_1$, $J$  disjoint linear chains  
$\Theta_i + G_i'$  in  $F_1$  of length  $g_i$ ($g_i \ge 1$) 
with  $\sum_{i=1}^J g_i = m$, 
and  $J$  disjoint sections  $C_i$  on  $Y_{\min}$  so that
$F_1 + C_i + \Theta_i + G_i'$  has the dual graph (2.3).
\newline
The surface $X$ is obtained by blowing up singular points (away
from  $\Theta_i$) and their infinitely near points on  $F_1$
(to get  $Y' \rightarrow Y_{\min}$), then blowing up a point  $q$ ($\notin C_i$) 
of a fibre  $F$ ($\ne F_1$)  on  $Y_{\min}$ (to get  $\pi' : X' \rightarrow Y'$) 
and finally blowing down smoothly the linear chains  $C_i + \Theta_i$.

\medskip\noindent
{\bf 2.9} {\bf Definition and Remark}. (1) $Y' \rightarrow Y_{\min}$
is not identical for otherwise  $h^0(-K_X) \ge h^0(-K_{X'}) = 1$.
Thus  $|M'|, P'$  in (2.4) are exactly the mobile, fixed part of
$|-2K_{X'} + 2E'|$.  $P'$  contains (but is contained in) the proper inverse transform 
(the total transform) of  $F_1$.

(2) Note that the unique member in  $|-K_{Y'}|$  contains  $\pi'_*(P')$
and also the support of the full fibre on  $Y'$  lying over  $F_1$
(cf. Lemma 1.9).  This and  $h^0(X', -K_{X'}) = 0$
explain why  $q \in F \ne F_1$.

(3) $X'$  is a Coble surface if and only if  $q$  is a singular
point of  $F$ (guaranteeing  $|-2K_{X'}| \ne \emptyset$).
If this is the case, then the  $X$  constructed as in Theorem 2.8 
with $m = 1$  is always a Coble surface (see 2.4 and Example 2.13);
it is also minimal if  $Y' \rightarrow Y_{\min}$  is a single blow-up
for then  $K_X^2 = -1$.  
See Example 2.10 for a situation where  $q$  has 2-dimensional
freedom to choose.

(4) A Coble surface  $X$  is of {\it Jacobian type}, or type (J) with respect to $E$ 
if
$X$  is equal to either  $X$  in Theorem 2.8
with an associated  $E$  there, or to
$X'$  with  $E' = E$  and  $q$  a singular point of  $F$.
In general, a Coble surface  $W$  is of {\it Jacobian type}
if there is a birational morphism  $W \rightarrow X$
such that  $X$  is of Jacobian type with respect to some  $E$;
$W$  is of elliptic type if it is either Jacobian or Halphen type (see 2.6);
$W$  is of {\it rational type} if it is not of elliptic type (cf. Example 2.11).
We can construct a minimal Coble surface which is of Halphen
type with respect to one curve but of Jacobian type with respect to another curve.

(5) The  $(-2)$-chain  $\Theta_i = \sum_{j=1}^{(g_i-1)} \Theta_i^{(j)}$  
meets only  $G_i'$  in  $F_1$, for otherwise  $M_1$, $G_i$  and one more
component  $P_1$  of  $P$  will share the same point,
which is absurd by Lemma 2.2.  Similarly, $G_i'\cdot \Theta_i = 1$.
In particular, $g_i \le 6$ ($= 6$ only when  $F_1$  is of type  $II^*$)
and  $g_i = 0$  when  $F_1$  is reduced.
When  $F_1$  is not reduced, 
it is impossible that  $g_i \ge 2$  for two  $i$  say  $i = 1,2$,
for otherwise the shortest chain in  $F_1$  linking  $G_1'$
and  $G_2'$  will give rise to a chain  $L$  in  $P$ 
(not a trivial fact; cf. Lemma 1.9 and (1), (2) above),
and hence to a loop  $M_1 + G_1 + L + G_2$  in  $M_1 + P$,
a contradiction to Lemma 2.2.

\par
Thus  $m = \sum_i g_i \le 6$.  Indeed, otherwise  $m \ge 7$,
$F_1$  is of type  $I_s$ ($s \ge m$) and  $C_i$
are sections intersecting different components of  $F_1$;
contracting all  $C_i$  and  $[m/2]$  components of  $F_1$,
we get a smooth rational surface  $V$  with  $K_V^2 = m + [m/2]
\ge 10$, a contradiction.
See Example 2.13 for the converse to Theorem 2.8.

\medskip\noindent
{\bf 2.10} {\bf Example}.  Let us give an example when $m = M^2 = 6$ occurs. 
Take two triples of non-concurrent lines
$(L_1,L_2,L_3), (L_4,L_5,L_6)$. Let us denote by $p_{ij}$ the intersection
point of the lines $L_i$ and $L_j$. We assume that   
$p_{12}\in L_4, p_{13} \in L_5,$ and 
$p_{23} \in L_6.$
The  curves $L_1+L_2+L_3$ and $L_4+L_5+L_6$ span a pencil of plane cubics with
nine base points 
$p_{16}, p_{25}, p_{34}, p_{12}, p_{13}, p_{23}$ and infinitely near points 
$p_{12}'\to p_{12}, p_{13}'\to p_{13}, p_{23}'\to p_{23}$
lying on the proper inverse transforms of the lines $L_4,L_5,L_6$. 
After blow up the base points we obtain
a Jacobian elliptic surface with  reducible fibres of type $I_6$ 
(its image in $\bbP^2$ is the union of
lines $L_1,L_2,L_3$) and  of type $I_3$ (its image in $\bbP^2$ is the union of
lines $L_4,L_5,L_6$). It has six disjoint sections corresponding to 
the six base-points $p_{16}, p_{25}, p_{34}, p_{12}', p_{13}', p_{23}'.$ 
If we blow down the six sections and 
blow up all 6 singular points of the fibre of type  $I_6$
(to get  $Y$), followed by the blow-up of 
a singular point  $q$  (to get a curve  $E$)  of the fibre of type  $I_3$,
we obtain a minimal Coble surface of 
Jacobian type with $M^2 = 6$  and  $K_X^2 = -1$.
Note that we can choose  $q$  to be any point as long as it is not on
the fibre of type  $I_6$ (to make sure that  $|-K_X| = \emptyset$), because
$\dim|-2K_X + 2E| = M^2 = 6 > 2$  always implies that  $|-2K_X| \ne \emptyset$.

\medskip \noindent
{\bf 2.11} {\bf Example}.  We construct a minimal
Coble surface  $X$  with two  disjoint  $(-1)$-curves  
$E_0$, $E_2$  such that  $X$  is of elliptic type with respect to the
first  $(-1)$-curve  $E_0$ (type(J) with  $m = g_1 = 2$) but
of rational type with respect to the second one  $E_2$
and fitting Case (2) with  $(m, k) = (0, 2)$  in Theorem 3.2.

Consider a minimal rational Jacobian surface $V$ with two fibres  $F_1, F_2$  of type 
$I_0^*$ (= $\tilde D_4$). One obtains this surface as the blow up of 9 
base points of the pencil of cubic curves spanned by the curve $L_1+L_2+L_3$, where 
$L_i$ are lines concurrent at a point $q$, and $H_1+2H_2$, where $H_1$ is a 
line through $q$ and $H_2$ is a line not containing $q$. It is easy to 
locate four disjoint sections  $E_i$  on $V$. Three of them come by blowing up 
infinitely near base points to the points in $H_2\cap L_i$, and the fourth one 
is blown up from an infinitely near base point to the point $q$. 
Let $C_i$ be the components of the fibre
$F_1$, intersecting $E_i$, and $C_i'$  the same for the other fibre  $F_2$. 
Let $X'$ be the blow-up of $V$ at two points lying on the 
multiple component of $F_1$ and one point $p$ lying on the multiple component 
of  $F_2$.  

Let $X$ be the blow-down of the section $E_1$ and the component $C_1$. 
This is a minimal Coble surface which is of elliptic type with respect 
to the exceptional curve  $E_0$  blown down to $p$. Now observe that 
if we blow down the section $E_2$ on  $X$  to get a surface $Y$, we can verify that
$|-2K_Y| = |2L|+ P$, where $|L|$ is the pencil of smooth rational curves 
linearly equivalent to the image on $Y$ of the component  $C_1'$.
Another member of $|L|$ entering as a component of an anti-bicanonical 
effective divisor is equal to the image of  $C_2+C_2'$. 
Thus $X$ is of rational type with respect to $E_2$.

\medskip\noindent
{\bf 2.12} {\bf Example.} Here we construct
examples of Coble surfaces  $X$  of Jacobian type
with respect to $E$  so that  $M^2 = 3$  in notation of Lemma 2.2.
Consider the union of three lines $L_i$  and a conic  $C$ (we may degenerate it into
the sum  $L_4 + L_5$  of two distinct lines) in  $\bbP^2$  such that
$C + \sum_{i=1}^3 L_i$  is of simple normal crossing.
Blowing up the 9 intersection points  $L_i \cap L_j, \, C \cap L_i$,
we obtain a surface  $Y$  with isolated  $|-K_Y|$  represented by the
proper inverse transform of  $\sum_{i=1}^3 L_i$.
Also we see that  $|-2K_Y|$  has the mobile
part defined by the linear system of cubics through 
the six intersection points  $C \cap L_i$. 
To be precise, $|-2K_Y| = |M| + \sum_i L_i', \, M \sim L + C'$,
where  $L$  is the pull-back of a general line
and  $L_i', C'$  the proper inverses of  $L_i, C$.
Now let  $X$  be the blow-up of  $Y$
at a point  $q$  not on the unique member of  $|-K_Y|$,
with  $E$  the exceptional curve. 
Then  $h^0(-2K_X) \ge h^0(-2K_Y) - 3 = 1$
and  $X$  is a minimal Coble surface of Jacobian type
with respect to $E$.

\medskip \noindent
{\bf Example 2.13.} We now give examples with  $g_i = 1$, 
kind of converse statement to Theorem 2.8 and Remark 2.9.  
The same idea can be applied to get examples with  $g_i \ge 2$
(see also Example 2.11).

Let  $Y_{\min}$  be a Jacobian minimal rational elliptic surface
with singular fibres  $F_1, F$.
Suppose that  $C_i$ ($1 \le i \le m$)  are  $m$  disjoint sections
meeting  $m$  different components  $G_i$  of  $F_1$.
We construct a blow-up  $Y' \rightarrow Y_{\min}$  in the following way:
it is the minimal blow-up of singular points and their infinitely
near points of  $F_1$  such that the proper inverses
of  $G_i$  all become  $(-4)$-curves on  $Y'$.
Let  $Y' \rightarrow Y$  be the blow-down of  $(-1)$-curves  $C_i$.

Then one can verify that  $|-K_Y|$  has only one member  $\sum_i G_i' + \Delta$,
where  $G_i'$  is the strict transform of  $G_i$  which is a  $(-3)$-curve
with  $G_i' \cdot \Delta = 2$,
where  $\Delta$  is effective and contractible to 
the divisor  $F_1 - \sum_i G_i$  and hence further to Du Val singularities
(to be precise, it is a set of a few smooth points when  $F_1$  is reduced).
We have also
$$|-2K_Y| = |M'| + P', \,\,\, P' = \sum_i G_i' + H',$$
where  $0 \le H' \le \Delta$,
where  $M'$  is the strict transform of a general full fibre
and hence a smooth elliptic curve with  $(M')^2 = m$,
Let  $\pi : X \rightarrow Y$
be the blow-up of a singular point  $q$  on the strict transform  $F'$  on  $Y$
of the second fibre  $F$.  Then  $X$  is a Coble surface
with  $M^2 = m$  in notation of Lemma 2.2, where  $M = \pi^*M'$.
Indeed, $|-K_X| = \emptyset$  for  $q$  is not on the unique member of 
$|-K_Y|$;
$|-2K_X| \ne \emptyset$  because  $F'$  is a member of  $|M'|$
with  $\roman{mult}_q F' \ge 2$ (Lemma 1.5).

\noindent
\head 3. The rational case \endhead

\noindent
{\bf 3.1} Now we shall consider the case of a Coble surface  $X$
and a  $(-1)$-curve  $E$  with  $p_a(-2K_X+2E) = 1$
(on a minimal Coble surface, any  $(-1)$-curve satisfies this),
such that the mobile part  $|M|$  of   
$|-2K_X+2E| = |M| + P$  satisfies  $p_a(M) =0$. 
As in Lemma 2.2, write  $M = kM_1$ ($k \ge 1$) with  $M_1 \cong \bbP^1$, 
and the fixed part as  $P = G + H = \sum_{i=1}^J g_iG_i + H$,
where  $G_i \ne G_j$  when  $i \ne j$.  
Write also  $H = \sum_i H_i$  where
$H_i = H_j$  is allowed.  We note that  $k \ge 2$  happens
only when  $M_1^2 = 0$.  Set  $m = M_1^2$. Let us state the 
theorem classifying all Coble surfaces of rational type.

\medskip\noindent
\plainproclaim 3.2 Theorem. There is a birational morphism  
$\tau : X \rightarrow Y_{\min}$ onto a minimal rational surface $Y_{\min}$, 
factoring as the blow down  $\pi : X \rightarrow Y$  of  $E$  and a morphism  
$\tau_y : Y \rightarrow Y_{\min}$, such that 
the direct image  $\Gamma := k{\overline M}_1 + \sum_{i=1}^J g_i {\overline G}_i +
{\overline H} \in |-2K_{Y_{\min}}|$  of  $kM_1 + G + H \in |-2K_X + 2E|$
are described as in one of the following Cases (1) - (16),
where to save notation, we use the same  $M_1, G_i, H_i$  to denote
their images  ${\overline M}_1, {\overline G}_i, {\overline H}_i$
on  $Y_{\min}$.

In Cases (1) - (9), $Y_{\min} = \bbP^2$  and hence
$\Gamma$  is a sextic.

\roster
\item $\Gamma = M_1 + 2G_1 + H_1$; $(m, k) = (0, 1)$;
$G_1$  is a conic, $M_1$  and  $H_1$  are distinct lines
meeting at  $p_1$; Supp $\Gamma$  is of simple normal crossing.

\item $\Gamma = kM_1 + 2G_1 + \sum_{i=1}^{4-k} H_i$;
$(m, k) = (0, 1), (0, 2)$; 
$M_1$  and  $H_i$  are lines through the same point  $p_1$;
$G_1$  is a line not through  $p_1$; $H_i = H_j$  is allowed
but  $M_1 \ne H_i$.

\item $\Gamma = M_1 + 2G_1 + H_1$, 
$(m, k) = (0, 1)$; $G_1$  is a conic; $M_1$  and  $H_1$  are distinct lines 
intersecting  $G_1$  transversally at the same point  $p_1$  and two other points.
\item $\Gamma = kM_1 + \sum_{i=1}^{6-k} H_i$;
$(m, k) = (0, k)$  with  ($1 \le k \le 6$);
$M_1, H_i$  are lines through the same point  $p_1$; 
$H_i = H_j$  is allowed but  $M_1 \ne H_i$.

\item $\Gamma = M_1 + \sum_{\ell=1}^J g_{\ell} G_{\ell}$;
$(m, k) = (1, 1)$;  $g_1 = 1, 2$; $2g_1 + \sum_{j=2}^J g_j = 5$;
$G_1$  is a conic; $M_1$  and  $G_j$ ($2 \le j \le J$)
are  $J$  distinct lines;
Sing $\sum_{\ell=1}^J G_{\ell}$  is disjoint from  $M_1$.

\item $\Gamma = M_1 + \sum_{i=1}^J g_i G_i$;
$(m, k) = (1, 1)$; $\sum_{i=1}^J g_i = 5$; 
$M_1$  and  $G_i$ ($1 \le i \le J$)
are  $J+1$  distinct lines;
$G_i$  and  $G_j$  share no common points on  $M_1$  when  $i \ne j$. 

\item $\Gamma = M_1 + 3G_1 + G_2$;
$(m, k) = (3, 1)$; $M_1$  is a conic;
$G_i$  are distinct lines; Supp $\Gamma$
is of simple normal crossing; let  $p_1$  be
a common point of  $M_1$  and  $G_2$.

\item $\Gamma = M_1 + \sum_{i=1}^J g_i G_i$;
$(m, k) = (3, 1)$; $g_1 = 1, 2$; $\sum_{i=1}^J g_i = 4$;
$M_1$  is a conic; $G_i$  are distinct lines;
all  $G_j$ ($2 \le j \le J$) intersect  $M_1$
transversally at the same point  $p_1$  and  $J-1$
other points; $G_1$  meets  $M_1$  at two distinct 
points not in  $M_1 \cap G_j$ ($j \ge 2$).

\item $\Gamma = M_1 + 4G_1$;
$(m, k) = (4, 1)$; $M_1$  is a conic; $G_1$
is a line intersecting  $M_1$  at two distinct points.

\item $Y_{\min} = \bbP^1 \times \bbP^1$;
$\Gamma = M_1 + \sum_{\ell=1}^J g_{\ell}G_{\ell}$;
$(m, k) = (2, 1)$; $g_1 = 1, 2$; $\sum_{i=2}^r g_i = \sum_{j=r+1}^J g_j 
= 3-g_1$; $M_1$  and  $G_1$  are sections (of both rulings)
of self intersection 2 and intersect each other
at two distinct points; $G_i$ ($2 \le i \le r$)
and  $G_j$ ($r+1 \le j \le J$)  are distinct
fibers of two different rulings such that 
Sing $\sum_{\ell=1}^J G_{\ell}$  is disjoint from  $M_1$.

\item $Y_{\min} = \bbP^1 \times \bbP^1$;
$\Gamma = M_1 + \sum_{\ell=1}^J g_{\ell} G_{\ell}$;
$(m, k) = (2, 1)$; $\sum_{i=1}^r g_i = \sum_{j=r+1}^J g_j = 3$;
$M_1$  is a section (of both rulings) of self intersection 2;
$G_i$ ($1 \le i \le r$) and  $G_j$ ($r+1 \le j \le J$)  are distinct 
fibers of two different rulings such that
$G_i \cap G_j \cap M_1 = \emptyset$.

\item $Y_{\min} = {\bold F}_2$;
$\Gamma = M_1 + \sum_{\ell=1}^J g_{\ell} G_{\ell} + hH_1$
($0 \le h \le 3$); $(m, k) = (2, 1)$;
$g_1 = 3 - h$; $\sum_{j=2}^J g_j = 2h$;
$H_1$  is the unique  $(-2)$-curve on  ${\bold F}_2$;
$M_1$  and  $G_1$  are two sections of self intersection $2$
and intersect each other at two distinct points;
$G_j$ ($2 \le j \le J$)  are distinct fibers not through  
$M_1 \cap G_1$; when  $h = 0$ (resp. $h = 3$), there is no such  $G_j$
(resp. no such  $G_1$).

\item $Y_{\min} = {\bold F}_b$ ($b \ge 2$); 
$\Gamma = kM_1 + 4G_1 + \sum_{i=1}^{2(b+2) - k} H_i$;
$(m, k) = (0, k)$  with  $1 \le k < 2(b+2)$;
$G_1$  is the unique  $(-b)$-curve;
$M_1$  and  $H_i$  are fibres; 
$H_i = H_j$  is allowed but  $H_i \ne M_1$.

\item $Y_{\min} = {\bold F}_{m-2}$;
$\Gamma = M_1 + 3G_1 + \sum_{j=2}^J g_j G_j$;
$(m, k) = (m, 1)$  with  $m \ge 3$; $\sum_{j=2}^J g_j = m+1$;
$G_j$ ($2 \le j \le J$)  are distinct fibres not through  $M_1 \cap G_1$;
$G_1$  is the negative section with  $G_1^2 = -(m-2)$;
$M_1$  is a section with  $M_1^2 = m$.

\item $Y_{\min} = {\bold F}_{m-4}$;
$\Gamma = M_1 + 3G_1 + \sum_{j=2}^J g_j G_j$;
$(m, k) = (m, 1)$  with  $m \ge 4$;
$\sum_{j=2}^J g_j = m-2$;
$G_j$ ($2 \le j \le J$) are distinct fibres 
of a fixed ruling not through  $M_1 \cap G_1$;
$G_1$  is a section with  $G_1^2 = -(m-4)$;
$M_1$  is a section with  $M_1^2 = m$
and meeting  $G_1$  at two distinct points.

\item $Y_{\min} = {\bold F}_m$; 
$\Gamma = M_1 + 3H_1 + \sum_{i=1}^J g_i G_i$;
$(m, k) = (m, 1)$  with  $m \ge 3$;
$\sum_{i=1}^J g_i = m+4$;
$H_1$  is the unique  $(-m)$-curve on  $Y_{\min}$;
$M_1$  is a section with  $M_1^2 = m$;
$G_i$ ($1 \le i \le J$) are distinct fibres.
\endroster

\medskip \noindent
{\bf 3.3 Remark} (1) Let  $\tau_0 : Y_0 \rightarrow Y_{\min}$  be
the blow-up of the point  $p_1$  in Cases (1)-(4), (7), (8);
and we set  $\tau_0 = \roman{id}$ for other cases.  Then  $\tau$
constructed in the proof factors through  $\tau_0$.
Moreover, $M_1$  on  $X$  is the total transform of the proper
inverse image on  $Y_0$  of  ${\overline M}_1$  on  $Y_{\min}$.
So the advantage of this classification is that
we can cook up a Coble surface by choosing
the right material: $(Y_{\min}, \Gamma)$
according to the customer's taste: like request
for  $M_1^2, \, h^0(-2K_X+2E), \, |kM_1|, \, G \cap M_1$, etc.

(2) $\sum_i \bar G_i$  in Case (5) or (6)
must be a non-reduced divisor (see Theorem 6.3 in \S 6).
In Cases (11) and (15) with  $m = 4$ (resp. Case (10)), 
$\tau_y : Y \rightarrow \bbP^1 \times \bbP^1$
factors through the blow-up of the intersection  $\bar G_i \cap \bar G_j$
of fibres of different rulings
for some  $i, j$, by the argument in the proof of Theorem 3.2
for Case (10) to deduce  $g_1 \le 2$
(resp. by the uniqueness of a loop, if exists, in  $M+P$
on  $X$; see Lemma 2.2).

(3) See Examples 2.11 and 4.8 and Remark 4.9 for the realizations 
of Case (2) with $(m, k) = (0, 2)$, and Cases (13) - (16).

\medskip
Let us start proving Theorem 3.2.
\medskip\noindent
{\bf 3.4}   There are two main cases to consider 
(Lemma 2.2):

\medskip
Case I: $|M| = |kM_1|$, where  $M_1 \cong \bbP^1$  and  $M_1^2 = 0$.

Case II: $M \cong \bbP^1$, $M^2 = m \ge 1$  and  $\dim |M| = m+1$ (Lemma 1.7).

\medskip\noindent
{\bf 3.5} We begin with Case I.
In notation of Lemma 2.2, we have  $G\cdot M_1 = \sum_{i=1}^J g_i G_i\cdot M_1 = 4$.
Since  $|M| + G + H = |-2K_X+2E|$  contains a divisor $2D$, 
where  $D \in |-K_X+E| =\pi^*|-K_Y|$
and the (multi-)sections  $G_i$  cannot be a component of a divisor from  $|M|$, 
we see that each  $g_i$  is even.  Thus either  
$G = 2G_1$, or  $G = 4G_1$,  or  $G=2G_1+2G_2$.

\par
Since  $E\cdot H = 0, E\cdot M_1 = 0$ (Lemma 2.3), $E, H$  are
all contained in fibres of the fibration given by  $|M_1|$.
Applying the blow-down  $\pi : X \rightarrow Y$, we
get  $|-2K_Y| \sim |\pi_*M| + \pi_*(G+H)$.
Intersecting this equality with a negative curve  $C$
on  $Y$, we see that either  $C \le \pi_*(G + H)$  or 
$C$  is a  $(-n)$-curve with  $n = 1,2$.
Let  $\tau : Y \rightarrow {\bold F}_b$  be a suitable 
smooth blow-down of  $(-1)$-curves in fibres of
the fibration given by  $|\pi_*M_1|$.  We will choose  $b$  later.

\medskip\noindent
{\bf 3.6} Suppose  $G = 2G_1$.  Then  $G_1$  is a double section
of the fibration given by  $|M_1|$.  Now  
$1 \ge p_a(kM_1 + G_1) \ge h^0((k-1)M_1) = k$  implies that  $k = 1$ 
(Lemma 1.2).  By 3.5 and the proof of Lemma 4.2 below, 
we can choose  $\tau$  so that  $b = 1$ (noting that  $K_Y^2 \le 0 < 8$).
Combining  $\pi, \tau$  and the blow-down  ${\bold F}_1 \rightarrow \bbP^2$,
we get a birational morphism  $X \rightarrow \bbP^2$.
Clearly, Case (1) occurs; indeed, deg $\Gamma = 6$  implies that
the image on  $\bbP^2$  of  $H$  is a line.

Suppose that  $G = 4G_1$.  Then  $G_1$  is a section 
of the fibration given by  $|M_1|$.  
This times, we can choose  $b = -\tau(G_1)^2 = -G_1^2 \ge 1$ (cf. Lemmas 1.8 and 2.2).
If  $b = 1$, combine  $\pi, \tau$  and the blow-down
${\bold F}_1 \rightarrow \bbP^2$  of  $\tau(G_1)$
and we get  $X \rightarrow \bbP^2$  fitting Case (4).
If  $b \ge 2$, then Case (13) occurs, and $k < 2(b+2)$  
by the reasonning as in 3.14.

Suppose that  $G = 2(G_1+G_2)$.  Then the  $G_i$ 
are sections of the fibration given by  $|M_1|$.
Since  $p_a(kM_1+G_1+G_2) \le 1$  and by Lemma 1.3, we have  
$(k, G_1\cdot G_2) = (1, 0), (2, 0), (1, 1)$.
By Lemma 1.8, we may choose  $b = -\tau(G_2)^2 = 1$
with  $\tau(G_1)\cdot \tau(G_2) = G_1\cdot G_2$.
Thus Case (2) or (3) occurs.

\medskip\noindent
{\bf 3.7} Next we consider Case II. We have $M = M_1$  and  
$\sum_i g_i G_i \cdot M = 4 + M^2$ (Lemma 2.2). 
Denote by  $\tau: X \to \bbP^{m+1}$  a morphism
given by the linear system  $|M|$ (cf. Lemma 1.7).
Since  $H \cdot M = 0 = E \cdot M = 0$,
the map  $\tau$  contracts  $H$  and factors through
the blow-down  $\pi : X \rightarrow Y$  of  $E$.

Since  $M^2 = m$, the image of  $\tau$  is a surface of degree  $\le m$ in
$\bbP^{m+1}$.  On the other hand, a non-degenerate surface in  $\bbP^{m+1}$ 
has degree $\ge m$.  So  $\tau$  is a birational morphism onto a surface 
$V$ of degree  $m$.  Such surfaces were classified by del Pezzo. 
According to his classification (see for example, [Re], p. 27),
$V$ is either  $\bbP^2$  ($m = 1$) or a Veronese surface  $V_4\subset \bbP^5$
($m = 4$), or a rational scroll $\bbF(a,n)$.
The latter surface is the image of a minimal ruled surface  $\bfF_n$  
under the map given by the linear system $|af+s_0|$, 
where  $f$  is the general fibre of the fixed ruling and  $s_0$  
a section with  $s_0^2  = -n$  and  
$m = 2a-n, a\ge n$.  Note that $\bbF(k,k)$ is the projective cone over 
a normal rational curve  $C\subset \bbP^k$  of degree  $k$. 
  
The following result follows easily from Lemma 2.2.

\medskip
\plainproclaim 3.8 Lemma.
\roster 
\item 
$G_i\cdot M \le 2$; if  $G_1\cdot M = 2$
then  $G_i\cdot M = 1$  and  $G_i \cap G_1 = \emptyset$  for all  $i \ge 2$.
\item
$G_i\cdot G_j \le 1$  for  $i \ne j$; if  $G_1\cdot G_2 = 1$
then  $G_i\cdot M = 1$  and  $G_i \cap (G_1 + G_2) = \emptyset$  for all  $i \ge 3$.
\endroster

\medskip
Now we shall treat possibilities of  $V$  in  3.7 one by one.

\medskip\noindent
{\bf 3.9}. Suppose that  $m = 4$  and  $V \subset \bbP^5$ 
is a Veronese surface.
Then  $M$  is the pull back of a conic in  $V$, viewed as a curve in  $\bbP^2$.
Hence  $M\cdot G_i \ge 2$  always holds.
This, together with Lemma 3.8, implies that  $G = 4G_1$,
the image  $\tau(G_1)$  is a line in  $V = \bbP^2$
and Case (9) occurs. 

\par 
Suppose that  $m = 1$.  Then  $V = \bbP^2$  and  $M$  is the pull back
of a line.  As in the case  $m = 4$, Lemma 3.8 implies that
Case (5) or (6) occurs.

\medskip\noindent 
{\bf 3.10}.  Now let us consider the remaining cases where $V = \bbF(a,n)$. 
First observe that in the case $V = \bbF(m, m), m \ge 2$, 
the map  $\tau : X \rightarrow V$  
factors through a birational morphism  ${\overline \tau}: X \to {\overline V}$, 
where  ${\bar V} = \bfF_m$. 
Now  $|M| + G + H$  is a subsystem of the pull back of the
linear system  $|mf+s_0|$.  

\par
So, in the case $V = \bbF(a, n), m \ge 2$, we have a map from  $X$  onto 
$\bfF_n$  such that  $|M| + G+H$  is a subsystem of the pull back of 
$|af+s_0|$, with  $2a-n = m$  and  $a \ge n$.
Let  ${\overline G}_i \sim a_if+b_is_0$  be the image of $G_i$, 
and  ${\overline H} = hs_0$  the image of  $H$, where  $h \ge 0$  and  $h \ge 1$
only when  $a = n = m$.
Since  $(2n+4)f+4s_0 \sim -2K_{\overline V}$  is linearly equivalent to 
the direct image of  $M + \sum_i g_i G_i + H$, we obtain
$$\sum_{i=1}^J g_ia_i+a = 2n+4, \,\,\,\, \sum_{i=1}^J g_ib_i + h = 3. \eqno (3.1)$$
In particular, $0 \le h \le 3$.

>From Lemma 3.8, we also obtain 
$G_i\cdot M = (a_if+b_is_0)\cdot (af+s_0) = a_i + (a-n) b_i = 1$ or  $2$.
Moreover, since  ${\overline G}_i$ is irreducible, $a_i\ge nb_i$, unless $(a_i,b_i) =
(0,1)$. This easily gives the following possible types:

\roster
\item $a_i = 1, b_i = 0; G_i\cdot M = 1$;
\item $a_i = 1, b_i = 1, n  = 0, a = 1, m = 2; G_i\cdot M = 2$;
\item $a_i = 2, b_i = 1, n  = 2, a = 2, m = 2; G_i\cdot M = 2$;
\item $a_i = 1, b_i = 1, n = 1, a = 2, m = 3; G_i\cdot M = 2$;
\item $a_i = 0, b_i = 1, a = n+1, m = n+2; G_i\cdot M = 1$;
\item $a_i = 0, b_i = 1, a = n+2, m = n+4; G_i\cdot M = 2$.
\endroster

On the other hand, by Lemma 2.2, $4+m = \sum_i g_i G_i\cdot M =
\sum_i g_i[a_i + (a-n) b_i]$.  Substituting (3.1) into this, we get
$$(a-n)[3 - \sum_{i=1}^J g_i b_i] = 0.  \eqno (3.2)$$
Hence either  $a = n = m$  and  $\sum_i g_i a_i = m+4$,  
or  $\sum_i g_i b_i = 3$  and  $h = 0$.

\noindent
Clearly, now we can divide into the following situations in 3.11-14.

\medskip\noindent
{\bf 3.11} For all  $1 \le i \le J$, Type 3.10 (1) occurs, i.e., $a_i = 1, b_i = 0$.
Then by (3.2) and (3.1), $\sum_i g_i = m+4, a = n = m, h = 3$.
This fits Case (12) ($m = 2$) or Case (16) ($m \ge 3$) of the theorem.  

In the following, we assume that for at least one  $i$,
Type 3.10 (1) does not occur.

\medskip\noindent
{\bf 3.12} $a = n = m$.  By 3.10, for each  $i$, either Type 3.10 (1) or (3)
occurs.  In view of Lemma 3.8, we may assume that  $G_1\cdot M_1 = 2$
(resp. $G_i\cdot M_i = 1$)  and Type 3.10 (3) (resp. (1)) occurs 
for  $i = 1$ (resp. for  $2 \le i \le J$).
Then  $n = 2, \, 2g_1 + \sum_{j \ge 2} g_j = 6, \, g_1 + h = 3$.
This is Case (12) with  $h \le 2$.

>From now on, we assume that  $a \ge n+1$  and hence  $h = 0$.

\medskip\noindent
{\bf 3.13} Suppose that for all  $i$, we have  $G_i\cdot M = 1$.
We may assume that
for  $1 \le i \le r; r \ge 1$ (resp. $r+1 \le j \le J$), Type 3.10 (5) 
(resp. (1)) occurs.  Thus  $a = n+1, m = n+2$, $\sum_{i=1}^r g_i = 3$,
$\sum_{j=r+1}^J g_j = m+1$.
If  $n = 0$, Case (11) occurs.
If  $n \ge 1$, then the uniqueness of a negative curve
on  ${\bold F}_n$  implies that  $r = 1$.  Hence  $g_1 = 3$  and
$\sum_{j\ge 2} g_j = m+1$.  Case (14) occurs.

\medskip\noindent
{\bf 3.14} In view of Lemma 3.8, we may assume now that  $G_1\cdot M = 2$, i.e.,
for  $i = 1$ Type 3.10 (2), (4) or (6) occurs,
and  $G_i\cdot M = 1$  for all  $i \ge 2$.  Then  $b_1 = 1$.
We may also assume that for  $2 \le i \le r$ ($r \ge 1$)
(resp. $r+1 \le i \le J$) Type 3.10 (5) (resp. (1)) occurs.
Thus  $g_1 + \sum_{i=2}^r g_i = 3$, $g_1 a_1 + \sum_{j=r+1}^J g_j + a = 2n+4$.

If Type 3.10 (6) occurs when  $i = 1$, then
$a = n+2, m = n+4$; hence  $r = 1$, $g_1 = 3$, $\sum_{j\ge 2} g_j = m-2$;
so Case (15) occurs.

If Type 3.10 (2) occurs when  $i = 1$, then  $a = 1$, 
$n = 0, m = 2, g_1 + \sum_{i=2}^r g_i = g_1 + \sum_{j\ge r+1} g_j = 3$.
So Case (10) occurs.  We note that  $g_1 \le 2$  for otherwise
$Y \rightarrow \bbP^1 \times \bbP^1$  is the blow-up of
points on  $G_1 \setminus M$  and their immediate infinitely near points
(cf. Lemma 1.9) and  $G_1 + M \in |-K_{Y_{\min}}|$  would give rise to 
a member in  $|-K_X|$ ($= \emptyset$), a contradiction. 

Suppose that Type 3.10 (4) occurs when  $i = 1$.  Then  $a = 2$,
$n = 1, m = 3, g_1 + \sum_{i=2}^r g_i = 3, g_1 + \sum_{j=r+1}^J g_j = 4$.
Thus, either  $r = 1$,
or  $r = 2$  and  ${\overline G}_2$  is the 
unique  $(-1)$-curve on  ${\bold F}_1$.
If  $r = 1$, then  $g_1 = 3, g_2 = 1$; we blow down the  $(-1)$-curve  $s_0$
and see that Case (7) occurs.
If  $r = 2$, then  $g_1+g_2 = 3, g_1 + \sum_{j \ge 3} g_j = 4$; 
we blow down the  $(-1)$-curve  ${\overline G_2}$
and see that Case (8) occurs.
This completes the proof of Theorem 3.2.

\head 4. Basic surfaces \endhead

\medskip\noindent
{\bf 4.1} A rational surface is called {\it basic} if it admits a 
birational morphism to $\bbP^2$ [Ha].  In the present section, 
we shall describe minimal Coble surfaces which are basic surfaces.

We start with the following well-known result:

\medskip
\plainproclaim 4.2  Lemma. Let $V$ be a rational surface. 
Suppose that $V$ does not have smooth rational curves with self-intersection $\le -3$. 
Then $V$ is a basic surface unless it is isomorphic to 
$\bfF_0$  or  $\bfF_2$.

{\sl Proof.} Let $\pi:V\to S$ be a birational morphism to a minimal rational surface  $S$. 
If $S\cong \bbP^2$ or $\bfF_1$, we are done. If $S \cong \bfF_b$ with $b \ge 3$, 
then the proper inverse transform of the negative section $s_0$  on 
$S$ is a curve on  $V$  with self-intersection $\le -b$, contradicting the assumption. 
If $S = \bfF_2$, then the same argument shows that $\pi$ is an isomorphism over the 
negative section  $s_0$.  So $\pi$ factors through a map $V'\to S$ 
which is the blow-up at a point on a fibre not lying on $s_0$. 
We blow down the proper transform of this fibre on $V'$ 
to get a morphism $V\to \bfF_{1}.$  The case  $S = \bfF_0$  is similar.

\plainproclaim 4.3 Theorem. Any Coble surface $X$  of 
elliptic type is basic.

{\sl Proof.} We may assume that  $X$  is of elliptic type with respect to some
$E$ (cf. Definition 2.9).  By Remark 2.6 and Theorem 2.8, $X$ is either one blow-up
of an Halphen surface of index 2, or is obtained by blowing up
a Jacobian elliptic surface  $Y_{\min}$  to get  $X'$  and then
blowing down linear chains of total length  $m \le 6$ (Remark 2.9).
An Halphen surface is a basic surface since it does not contain
smooth rational curves with self-intersection $\le -3$ (this immediately
follows from the formula for the canonical class). 
So, in the Halphen case  $X$  is basic.
For the Jacobian case, note that the exceptional divisors of  
$X' \rightarrow Y_{\min}$  and  $X' \rightarrow X$  are disjoint.
Thus there are smooth blow-downs  $Y_{\min} \rightarrow Z$
and  $Y \rightarrow Z$  fitting the following commutative diagram
(the rectangular part):
$$\hskip 1.3pc X' \overset{\pi'}\to\longrightarrow Y' \longrightarrow Y_{\min}$$
$$\downarrow \hskip 2.5pc \downarrow \hskip 2.5pc \downarrow$$
$$\hskip 3.3pc X \overset{\pi}\to\longrightarrow Y 
\hskip 0.2pc \longrightarrow \hskip 0.2pc Z \hskip 0.5pc \longrightarrow \bbP^2.$$
Since  $Y_{\min}$  satisfies the hypothesis of Lemma 4.2 so does  $Z$.
Therefore, there is a smooth blow-down  $Z \rightarrow \bbP^2$
because  $K_Z^2 = m + K_{Y_{\min}}^2 < 8$.  Theorem 4.3 follows.

\medskip\noindent
{\bf 4.4} Now let us assume that $X$ is a Coble surface of rational 
type with respect to a  $(-1)$-curve  $E$.  Write  $|-2K_X+2E| = |M|+P$ as in 2.1.
By Theorem 3.2 and Remark 3.3, $X$  is basic unless one of Cases (12)-(16)
occurs.  In these five cases we have a birational morphism  
$X \rightarrow {\bold F}_d$.  If Case (13) occurs, then  $M = kM_1$
and  $|M_1|$  is a free pencil of rational curves.
There is no upper bound for  $k$ (see Example 4.8 below);
of course if  $X$  dominates  $\bbP^2$,
via the blow-down  $\pi: X \rightarrow Y$  of  $E$,
then  $k \le 6$.  If one of Cases (12), (14)-(16) in Theorem 3.2 occurs,
then  $M$  is a smooth rational curve with  $m = M^2 \ge 2$;
again there is no upper bound for  $m$ (see Remark 4.9 below);
however if  $X$  dominates  $\bbP^2$, via  $\pi$,
then clearly  $M^2 \le 36$.  We can do much better.  We shall start with the following:

\plainproclaim 4.5 Lemma. Let $C$ be an irreducible rational plane
curve of degree $4
\le d \le 6$. Assume that  $C$ does not have a singular point of
multiplicity $d-1$ and, in the case $d = 6$, there is a a point of
multiplicity $\ge 3$.  Then there exists a Cremona transformation with 
fundamental points among singular points of $C$ such that the image of
$C$ is a curve of degree $\le 3$.

{\sl Proof.} Let $m_1\ge m_2\ge \ldots\ge m_k$ be the multiplicities of singular points
of $C$ (including infinitely near points). Consider the
vector $(d;m_1,\ldots,m_k)$. 

Case $d = 4$. The possible multiplicities of singular points are 
$(m_1,\ldots,m_k) = (4;2,2,2)$. We apply the standard quadratic 
Cremona transformation $T$ with centers at the singular points to get a conic.

Case $d = 5$. Then $(d;m_1,\ldots,m_k) = 
(5;3,2,2,2);(5;2,2,2,2,2,2)$.  In the first case, applying $T$ as above with centers at
the first three points, we get
$(d';m_1',\ldots,m_k') = (3;2)$. In the second case, we use the Cremona
transformation given by the linear system of quintics through the
singular points of the curve. We get a line.

Case $d = 6$. Assume 
$C$ has a point of multiplicity $4$. Then 
$(d;m_1,\ldots,m_k) = (6; 4,2,2,2,2)$.  We make a standard
Cremona transformation at the first three points.
Then $(d';m_1',\ldots,m_k')$ $=$ $(4;2,2,2)$.  Applying again the
standard quadratic Cremona transformation we get a conic.
 
Assume $C$ has a point of
multiplicity $3$ but no points of multiplicity $4$. Then 
$(d;m_1,\ldots,m_k) = (6;3,3,3,2), (6;3,3,2,2,2,2), (6;3,2,2,2,2,2,2,2).$
Again we make the standard quadratic Cremona transformation at the first three
points. We get $(d';m_1',\ldots,m_k') = (3;2), (4;2,2,2),
(5;2,2,2,2,2,2).$ In the second case we apply again the standard qyadratic Cremona
transformation to get a nonsingular conic. In the third case we apply the Cremona
transformation given by quintics through the six singular points of the
curve. We get a line. 

\plainproclaim 4.6 Proposition. Assume that a Coble surface  $X$  admits
a birational morphism  $\tau: X \rightarrow \bbP^2$  with  $E$  a $(-1)$-curve on 
it blown down.  Suppose further that  $X$  is of rational type
with respect to  $E$, and write  $|-2K_X+2E| = |M| + P$  as in 2.1. 
Then  $M^2 \le 5$ (the equality is realizable).

{\sl Proof.} We may assume that $M^2 \ge 1$ and hence the general member 
$M$  of  $|M|$  is a smooth rational curve.
Let  $\pi : X \rightarrow Y$  be the blow-down
of  $E$, which is an isomorphism in a neighbourhood of
the divisor  $M + P$  which is disjoint from  $E$ (Lemma 2.3).  
Then  $\tau$  is the composition
of  $\pi$  and a birational morphism  $\tau_Y : Y \rightarrow {\bbP^2}$.
We have  $|-2K_Y| = |\pi(M)| + \pi(P)$ (Lemma 2.3).

First observe that  ${\overline M} := \tau(M) = \tau_Y(\pi(M))$  is a component 
of the sextic  $D := \tau_*(M+P) \in |-2K_{\bbP^2}|$,
and  $\tau : M \rightarrow \bar M$  is a resolution of the rational curve  $\bar M$.
In particular, $d := \deg \bar M \le 6$.
If  $\bar M$  has at worst  $r$  double singular points
(this is true when  $d \le 3$),
then  $r = (d-1)(d-2)/2$  and  $M^2 \le d^2 - 4r \le 5$.

Thus, we may assume that  $d = 4, 5, 6$  and  $\bar M$
has a singular point of multiplicity $\ge 3$.
Note that the surface  $Y$  is obtained 
from $\bbP^2$ by successive blow-ups of singular points of
effective anti-bicanonical divisors. 
Let  $T: \bbP^2 \to \bbP^2$  be a Cremona transformation with fundamental points 
in the set  $\Sigma$  of indeterminancies of the rational map
$\tau_Y^{-1}$.  Clearly, Sing $D \subseteq \Sigma$.
Composing  $T$  with $\tau_Y$, we get another birational morphism  
$Y \to \bbP^2$  such that  $D$  is replaced with
the image of  $D$  under $T$.  

Let us show that this could be used to reduce our proof to
the case when  $d = \deg \bar M \le 3$.
Consider first the case  $d = 6$.  
If  $D = \bar M$  does not have a point of multiplicity  $5$, we apply Lemma 4.5 
to get a Cremona transformation  $T$  such that the image of  $D$  is a cubic. 
If  $D$  has a point $p$ of multiplicity 5, 
then the indeterminacy set  $\Sigma$  of  $\tau_Y^{-1}$  consists of
$p$  and its infinitely near points.
Applying Lemma 1.9 repeatedly, we see that  $|-K_Y|$
contains a member  $3F + 2E_0 +$ (an effective divisor),
where  $E_0$  is the proper inverse transform of the exceptional
curve lying over  $p$  and  $F$  is a smooth fibre
of a  $\bbP^1$-fibration whose image on  $\bbP^2$  
is a line through  $p$.  This implies  $|-K_X| \ne \emptyset$,
a contradiction.

Consider the case  $d = 5$.  Then the residual component of $\bar M$ in
$D$ is a line  $L$.  If all singular points of  $\bar M$  are of
multiplicity  $\le 3$, then applying the previous lemma, we reduce $\bar M$ to a
curve of degree $\le 3$. 
If  $\bar M$ has a point  $p$  of multiplicity 4, 
we apply the standard Cremona transformation with fundamental points at $p$
and two points from the set  $\Sigma \cap L\cap \bar M$; 
for the existence of these two points, we note that the proper inverse transform on  $Y$  of
$L$  and  $M$  should meet each other at most twice (Lemma 2.2), while
$L \cdot M = 5$, whence we need to blow up at least
3 points in  $L \cap M$ (including infinitely near).
This will transform  $\bar M$  to a quartic with a triple point. 

If $\bar M$ is a quartic with a triple point  $p$, then the
residual curve in  $D$  is a conic $Q$ (possibly a double line  $2L$). 
We apply the standard Cremona transformation with fundamental points at  $p$
and two points from the set  $\Sigma \cap Q \cap \bar M$, 
which exist by the above reasonning.
This will transform $\bar M$ to a cubic. 
Thus we have reduced to the case  $d \le 3$
and Proposition 4.6 is proved.

\medskip
\plainproclaim 4.7 Theorem. Any Coble surface  $X$  with
$h^0(-2K_X) \ge 7$  is not basic (see Example 4.8 below
for  $X$  with arbitrarily large anti-bicanonical dimension).

{\sl Proof.} Let  $W$  be a Coble surface with  $h^0(-2K_W) \ge 7$.
Suppose the contrary that there is a birational morphism  
$W \rightarrow \bbP^2$.
Clearly this map factors through  $W \rightarrow X$  with
$X$  minimal Coble.  Note also that  $h^0(-2K_X) \ge h^0(-2K_W) \ge 7$.
Take any  $(-1)$-curve  $E$  on  $X$  blown down by the map
$X \rightarrow \bbP^2$.  

Write  $|-2K_X + 2E| = |M| + P$  as in 2.1.  If  $p_a(M) = 1$
or  $p_a(M) = 0$  with  $M = kM_1$  and  $M^2 = 0$, 
then  $h^0(X, -2K_X) \le h^0(X, M) - 1 = M^2 \le 6$, 
or  $h^0(X, -2K_X) \le h^0(X, kM_1) - 1 = k \le \deg (-2K_{\bbP^2}) = 6$
(cf. Lemma 1.7, Remark 2.9), a contradiction,
where we used the fact that  $E$  is not in the fixed part
of  $|-2K_X + 2E|$ (Lemma 2.3).

Therefore, the hypothesis of Proposition 4.6 is satisfied
and we have  $h^0(-2K_X) \le h^0(X, M) - 1 = M^2 + 1 \le 6$.
This contradiction proves Theorem 4.7.

\medskip
\plainproclaim 4.8 Proposition. Given any integers 
$a, b$ with  $a \ge 4$  and  $b \ge 2a$,
there is a Coble surface  $X$
which does not admit any birational morphism
$X \rightarrow \bbF_d$, where   $d \le  a - 3$
and which satisfies  $h^0(X, -2K_X) = a$  and  $K_X^2 = 4 - a - b$.

We prove this result by constructing examples fitting Case (13)
of Theorem 3.2.

\medskip\noindent
{\bf 4.9 Example}. Let  $n, b, t$  be non-negative integers satisfying:
$n \ge 3, n \ge t, b \ge t + 2(n-1)$.  

Let $s_0$ be the negative section on  $S = \bfF_b$ 
with self-intersection $-b$ and $F$  a fibre.
Take distinct fibers  $F_\ell$ ($1 \le \ell \le b - t - (n-2)$).
Then we can write
$$-K_{S} = 2s_0 + \sum_{i=1}^r F_i + 2 \sum_{j=r+1}^{r+s} F_j +
3 \sum_{k=r+s+1}^{r+s+t} F_k,$$
$$-2K_{S} = nF + 4s_0 + 2 \sum_{i=1}^r F_i + 
3 \sum_{j=r+1}^{r+s} F_j +
5 \sum_{k=r+s+1}^{r+s+t} F_k,$$
where we set
$$r = b - t - 2(n-1), \,\,\, s = n - t.$$
Let  $\sigma : Y \rightarrow S$  be the composite of 
the blow-ups of smooth points
on  $F_\ell \setminus s_0$  and their infinitely near points
such that  
$$\sigma^*(F_i) = H_i + J_i, \,\,
\sigma^*(F_j) - J_j = H_j + 2E_j + 2B_j,$$
$$\sigma^*(F_k) - J_k = H_k + 2E_k + 2B_k + 2C_k + 2D_k$$
have the following dual graphs:
$$(-1) - (-1),$$
$$(-2) - (-2) - (-1),$$
$$(-2) - (-2) - (-2) - (-2) - (-1).$$
Here  $H_{\ell}$ is the proper inverse transform of $F_\ell$, and  $J_u$ 
($u = j$ or  $u = k$) is a $(-2)$-curve with  $J_u\cdot E_u = 1$.

Then  $-K_Y$  equals
$$2G_1 + \sum_i H_i + \sum_j (2H_j + 2E_j + J_j + B_j) + 
\sum_k (3H_k + 4E_k + 2J_k + 3B_k + 2C_k + D_k),$$
and  $-2K_Y = nM_1 + 4G_1 + H$  where  $G_1 := \sigma^*(s_0)$, 
$M_1 := \sigma^*(F)$. Here we set
$$H = 2\sum_i H_i + \sum_j (3H_j + 2E_j + J_j) + 
\sum_k (5H_k + 6E_k + 3J_k + 4B_k + 2C_k).$$

Let $q \in Y$ be a point which either lies on $F_0 \setminus G_1$
with a smooth fibre  $F_0 \sim M_1$,
or on $J_i \setminus H_i, i \le r$. Let $\pi: X \to Y$ be the blow-up 
of $Y$ at $q$ and $E$ the exceptional curve.  We claim:

\roster
\item $X$ is a Coble surface.

\item $|-2K_X+2E| = |M|+G$, where 
$M = \pi^*(nM_1) = n(\sigma \circ \pi)^*(F)$  and 
$G = 4\pi^*G_1+ \pi^*H$.

\item
$h^0(X, M) = n+1$  and  $h^0(X, -2K_X) = n-1$.

\item $K_X^2 = 5-(n+t+b)$.
\newline
In the following we assume that  $n+t\ge 5$  and  $q$  lies on a smooth fibre  $F_0$.
 
\item $X$ does not admit a birational morphism 
to $\bfF_d$ with $d \le n+t-4$. In particular, $X$ is not a basic surface.
Moreover, all negative curves ($\ne \pi^{-1}G_1$) are contained in fibres.

\item
Suppose that  $r = 0$, i.e., $b = t + 2(n-1)$.  
Let  $X \rightarrow X_{\min}$  be the blow-down of
the proper inverse transform of  $F_0$.  Then  $X_{\min}$  is a 
minimal Coble surface.  

\item
Denote by  $B_j'$ (when  $s > 0$), $D_k'$ (when  $t > 0$)  
the proper images on  $X_{\min}$  of  $B_j, D_k$.  Then both mobile parts of
$|-2K_{X_{\min}}+2B_j'|$  and  $|-2K_{X_{\min}}+2D_k'|$  are equal to
$|(n+1)F|$  with  $F$  denoting a full fibre of the induced
$\bbP^1$-fibration on  $X_{\min}$.
\endroster

For (2), we only need to show that  $|-2K_Y| = |nM_1| + 4G_1 + H$.
To do so, we use the fact that for an effective divisor  $L$
and irreducible divisors  $N_i$, if  $L \cdot N_1 < 0, (L - N_1) \cdot N_2 < 0,
\cdots, (L - \sum_{i=1}^{v-1} N_i) \cdot N_v < 0$, then
$\sum N_i$  is a partial fixed part of  $|L|$.
Inductively, one can verify that  $|-2K_Y|$  contains the following
as its partial fixed part:
$$G_1 + \sum (H_u + E_u + J_u + B_k + C_k) + G_1 + \sum (H_u + E_u) + G_1 +$$
$$\sum (H_i + H_u + B_k + E_k + J_k + E_k + H_k) + G_1 + \sum H_i,$$
which is equal to  $4G_1 + 2 \sum H_i +$ (other components).
Since  $-2K_Y - (4G_1 + 2 \sum_i H_i)$  is a disjoint union of
$nM_1$  and a negative definite divisor contained in fibers,
$|nM_1|$  is the mobile part of  $|-2K_Y|$.  

For the last part of (5), we assume  $n \ge 5$  for simplicity.
Note that  $-2K_X \sim (n-2)\pi^*M_1 + G + 2 \pi^{-1}(F_0)$.
Suppose the contrary that  $C$ ($\ne \pi^{-1}G_1$)  is a negative curve
not contained in fibres.  Then  $C \cdot (-2K_X) \ge (n-2) M_1 \cdot \pi_*C
\ge n-2 \ge 3$.  This leads to that 
$2p_a(C) - 2 = C^2 + C \cdot K_X \le -1 - 2$, a contradiction.
The rest of the Claim can now be verified with patience.

\medskip\noindent
{\bf 4.10 Remark}. For each  $N = 14, 15, 16$,
we can construct similar Coble surfaces  $X$  fitting
Case (N) of Theorem 3.2 (as well as minimal Coble surface
$X_{\min}$  obtained as a single blow-down of  $X$) 
and with arbitrarily large  $-K_X^2$  and
$h^0(-2K_X)$  but with no birational morphism  $X \rightarrow \bbP^2$
(cf. 4.7 and 5.6).

\head 5. Coble sextics \endhead

\medskip
Let $X$ be a basic Coble surface.  So there is a birational morphism  $X \to \bbP^2$.
The image of any divisor $D \in |-2K_X|$ in $\bbP^2$ is a member of
$|-2K_{\bbP^2}|$, whence a plane sextic. A plane sextic which is the image of
an anti-bicanonical divisor of a basic Coble surface  (which we may assume
minimal) will be called a {\it Coble sextic}. In this section we shall
describe Coble sextics.

\medskip\noindent
{\bf 5.1} Assume that $X$ is a minimal Coble surface of Halphen type.  Then it is
obtained by blowing up a singular point of a fibre on an Halphen surface $V$
of index 2. We have already explained that $V$ is a basic surface. The image
on  $\bbP^2$  of the pencil of elliptic curves on  $V$  is an {\it Halphen pencil} of
index 2 of elliptic curves of degree 6 with 9 double base points, including
infinitely near. There is a unique plane cubic $C$ through the base points,
and the base points add up to a non-trivial 2-torsion point on the cubic with
one of the inflection points chosen as the origin (see [CD,Do1]). Even when the
cubic is a nodal curve, this makes sense.  The cubic $C$ taken with
multiplicity 2 is a member of the Halphen pencil. 

A Coble sextic corresponding to
$X$ is a member of an Halphen pencil of index 2  which has a singular point
$p$  such that a preimage on  $X$  of  $p$  is also a singular
point of a fibre dominating a member ($\ne 2C$) of the pencil on  $\bbP^2$.
The classical Coble sextic is of this type. The Halphen pencil has 9 distinct
double base point, and an irreducible  member of the pencil with an extra
singular point  $p$  is a rational sextic with 10 nodes. 

\medskip\noindent
{\bf 5.2} Assume now that $X$ is a minimal Coble surface of Jacobian type
as described in Theorem 2.8.  We use the commutative diagram in
the proof of Theorem 4.3.  

Suppose that the center  $q$  of the blow-up  $\pi' : X' \rightarrow Y'$  
is a singular point of the fibre  $F$  (this is not always true as shown in Example 2.10).
Then by Remark 2.9, $X'$  is a Coble surface and a member  $D'$  of  $|-2K_{X'}|$
is of the form  $D' := F_1' + F' +$ (an effective divisor contractible
by the map  $X' \rightarrow Y_{\min}$), where  $F_1', F'$  are
the proper inverse transforms of the two distinct fibres  $F_1, F$  on  $Y_{\min}$.
The image  $D$  on  $X$
of this  $D'$  is a member of  $|-2K_X|$.  Now the commutatitve
diagram in Theorem 4.3 shows that the image of  $D$  under the map
$X \rightarrow \bbP^2$, is equal to the image of  $F_1 + F$  under the
map  $Y_{\min} \rightarrow \bbP^2$.
Thus  $X$  or rather its anti-bicanonical divisor  $D$,
defines a Coble curve which is the union of two singular members of
a cubic pencil dominated by the elliptic fibration on  $Y_{\min}$.

For general  $q$  in  $F$, as above, the sextic image
$\Sigma$  of  $M + P \in |-2K_X + 2E|$ (or  equivalently of  $\pi_*(M+P) \in |-2K_Y|$)
under the birational morphism  $X \overset\pi \to\rightarrow Y \rightarrow \bbP^2$,  
is equal to the image of  $M_1' + P'$ (see (2.4)) under the map
$X' \rightarrow \bbP^2$, and hence equal to the image of  $F_1+F$
under the map  $Y_{\min} \rightarrow \bbP^2$ (cf. Remark 2.9).
If either  $F$  is smooth elliptic or both  $F$  and  $F_1$
have irreducible images on  $\bbP^2$, such  $\Sigma$  would
never be realized from a Coble surface of rational type
(cf. Corollary 5.5 below).

\medskip\noindent
{\bf 5.3 Remark.} The birational morphism $X\to \bbP^2$ constructed in Theorem 4.3
is not unique as the following example shows. Let $C_5$ be a plane curve of
degree 5 with six nodes. Let $L$ be a line intersecting $C_5$ at five
distinct points. Let us show that the sextic $C_5+L$ is a Coble sextic
obtained from a Coble surface of Jacobian type. Let $f:Y'\to
\bbP^2$ be the blow-up of 5 nodes $p_i$ of $C_5$ and four common points $q_j$
of  $C_5$ and $L$. 
The surface $Y'$ has an elliptic pencil  $\Lambda$  
spanned by the proper transform $F_1$ of 
$C_5$ and the union $F_2 = L' + 2C_2'$ of the proper transforms 
of $L$ and the double conic  $2C_2$ through the  points $p_i$. 
The pre-image of the point $q\in C_5\cap L$  different from
$q_j$'s is the unique base point of the pencil  $\Lambda$. The curves $F_1$ and $F_2$
are  singular members of the pencil. The singular point of $F_1$ is the
pre-image of the node $p$ of $C_5$ different from the points $p_i$'s. 
The $C_2'$ is now a $(-1)$-curve. 
Let $\gamma: X\to Y'$ be the blow-up of the singular point of $F_1$  with  $E_1$
the exceptional curve. It is easy to see that $X$ is a minimal
Coble surface of Jacobian type with respect to $E_1$  and also  $C_2'$
and with  $|-2K_X| = \{C_5' + L'\}$, where  $C_5'$
is the proper inverse of  $C_5$ (or  $F_1$).
The image of the anti-bicanonical divisor of  $X$ in
$\bbP^2$, under the map $f \circ \gamma : X \to \bbP^2$, is equal to $D_6 = C_5+L$
with 11 nodes $p_i, i = 1,\ldots, 5, q_j , j = 1,\ldots, 4$ and $p, q$. 

On the other hand, following Theorem 4.3, we blow down  $E_1$,
the  $(-1)$-curve  $C_2'$  in  $F_2$ (to get  $Y_{\min}$  
after further blow up the base point of  $\Lambda$) and also
sections and fibre components on  $Y_{\min}$, we get a new
birational morphism  $\sigma : X \rightarrow \bbP^2$,
which maps the anti-bicanonical divisor of  $X$  onto the union
of two nodal cubics (the images on this ``new''  $\bbP^2$  of  $C_5', L'$).
The two different Coble sextics on two ``different''  $\bbP^2$  derived from 
the same surface $X$ are related by the Cremona transformation of $\bbP^2$ 
defined by the two different birational morphisms from
$X$  to  $\bbP^2$. It can be given by the linear system of 
quintics with double points at $q,p_i, i = 1,\ldots, 5$,
if one chooses  $\sigma$  properly.

\bigskip
Next we consider Coble surface  $X$  of rational type with respect to a
$(-1)$-curve  $E$  on it.  As in Lemma 2.2, write  $|-2K_X+2E| = |M| + P$, $M = kM_1$,
$M_1 \cong \bbP^1$, $P = G + H$, $G = \sum_{i=1}^J g_i G_i$, $H = \sum_j H_j$.
We note that if  $\sigma : X \rightarrow \bbP^2$  is a birational morphism,
factoring as the blow-down  $\pi : X \rightarrow Y$  of the curve  $E$
and a morphism  $\sigma_y : Y \rightarrow \bbP^2$, then  $\Sigma := \sigma_*(M+P)$
is a sextic plane curve and equal to the  $\sigma_y$-image of the member
$\pi_*(M+P)$  in  $|-2K_Y|$.  We shall prove:

\medskip \noindent  
\plainproclaim Theorem 5.4. Assume that  $X$  is a basic surface of
rational type with  $E$  blown down by the map onto  ${\bbP^2}$.
Then there is a (possibly new) birational morphism  $\sigma: X \rightarrow
{\bbP^2}$  with  $E$  also blown down by it, such that  $\sigma$  and the 
sextic  $\Sigma = \sigma_*(kM_1 + G + H)$  are equal to one of the following,
where for simplicity, we employ the same symbols  $M_1, G_i, H_i$
to denote their  $\sigma$-images  ${\hat M}_1, {\hat G}_i, {\hat H}_i$
in  ${\bbP^2}$:
 
\roster
\item
$\sigma$  and  $\Sigma$  are identical to  $\tau$, $\Gamma$  in
one of Cases (1)-(9) in Theorem 3.2;
so  $\Sigma$  is a union of lines and conics.

\item
$X, E$  fit Case (13) of Theorem 3.2 with  $1 \le k \le 6$;
$\sigma$  is the blow-down  $X \rightarrow {\bold F}_1$  of  $E$
and all curves in fibres of the  $\bbP^1$-fibration given by  
$|M_1|$  so that  $G_1$  becomes the  $(-1)$-curve on  $\bold F_1$,
followed by the blow-down  $\bold F_1 \rightarrow {\bbP^2}$
of  $G_1$; 
$\Sigma = k{\hat M}_1 + \sum_{j=1}^{6-k} {\hat H}_j$, 
where  ${\hat M}_1$, ${\hat H}_j$  are concurrent lines,
with  ${\hat M}_1 \ne {\hat H}_j$, but  ${\hat H}_i = {\hat H}_j$
possible; the $H_j$  here may be different from the  $H_i$  in Theorem 3.2.

\item
$X, E$  fit Case (14) of Theorem 3.2 with  $m = 3$; $\sigma$
is the composition of  $\tau : X \rightarrow \bold F_1$
and the blow-down  $\bold F_1 \rightarrow \bbP^2$  of  $\bar G_1 = \tau(G_1)$;
$\Sigma = \hat M + \sum_{j=2}^J g_j \hat G_j$  with  
$\sum_{j=2}^J g_j = 4$,
where  ${\hat G_j}$  are lines concurrent at  $p$
and  $\hat M$  is a conic through  $p$  and transversal to all  $\hat G_j$.
 
\item
$X, E$  fit Case (15) of Theorem 3.2 with  $m = 5$; $\sigma$
is the composition of  $\tau : X \rightarrow \bold F_1$
and the blow-down  $\bold F_1 \rightarrow \bbP^2$  of  $\bar G_1$;
$\Sigma = \hat M + \sum_{j=2}^J g_j \hat G_j$  with  
$\sum_{j=2}^J g_j = 3$,
where  $\hat M$  a cubic with a node at  $p$, where
${\hat G_j}$  are lines through  $p$  and transversal to both
tangents of  $\hat M$  at  $p$.

\item
The well-defined morphism  $\sigma$  is the composition of  
$\tau : X \rightarrow {\bbP^2}$  in Case (10), or (11), or (15) 
with  $m = 4$  of Theorem 3.2, the blow-up of an intersection point
${\bar G}_i \cap {\bar G}_j$  of fibres of two different rulings
with exceptional divisor  $D$  and the blow-down of the proper inverses of 
these two fibres (Remark 3.3);
$\Sigma = {\hat \Gamma} + (g_i+g_j - 2){\hat D}$, where  ${\hat \Gamma}$
is the strict transform of  $\Gamma$  and  ${\hat D}$  the image of  $D$;
so  $\Sigma$  is a union of lines, conics and at most one nodal
cubic (only in Case (15), and then the  $\Sigma$  here is the same
as the one in (4) above).

\item
$X, E$  fit Case (12)  of Theorem 3.2; 
there is a section  $C$  of the
$\bbP^1$-fibration on  $X$  induced from the one on  ${\bold F}_2$,
with  $\pi(C)^2 = -1$, $C \cap (G+H) = \emptyset$  and  $C \cdot M = 2$;
$\sigma$  is the blow-down  $X \rightarrow {\bold F}_1$  of  $E$
and all curves in fibres disjoint from  $C$  followed by the
blow-down  $\bold F_1 \rightarrow {\bbP^2}$  of  $C$; 
$\Sigma = {\hat M} + (3-h) {\hat G}_1 + h {\hat H}_1$, 
where  $0 \le h \le 3$, ${\hat M}$  is a cubic with a node at  $p$  
and  ${\hat G}_1$, ${\hat H}_1$  are distinct lines not through  $p$.
 
\item
$X, E$  fit Case (13) of Theorem 3.2 with  $k = 1, 2$;
there is a section  $C$  of the  $\bbP^1$-fibration on
$X$  given by  $|M_1|$, with  $\pi(C)^2 = -1$, $C \cap G = \emptyset$
and  $C \cdot H = (2-k)$;
$\sigma$  is the blow-down  $X \rightarrow {\bold F}_1$  of  $E$
and all curves in fibres disjoint from  $C$  followed by the
blow-down  $\bold F_1 \rightarrow {\bbP^2}$  of  $C$; 
$\Sigma = k{\hat M}_1 + 4 {\hat G}_1 + (2-k) {\hat H}_1$, 
where  ${\hat M}_1$, ${\hat G}_1$, ${\hat H}_1$  are non-concurrent lines;
the  $H_1$  here may be different from any  $H_i$  in Theorem 3.2.
\endroster

\medskip \noindent
\plainproclaim 5.5 Corollary. With the assumptions in Theorem 5.1, we have:
\roster
\item 
The plane sextic  $\Sigma$  is a union of lines, conics and at most one nodal cubic;
moreover, if a cubic does appear in  $\Sigma$  then it is the image
of the mobile part  $M$  of  $|-2K_X + 2E|$.
\item
$M^2 \le 5$  holds; if  $M^2 = 5$  then (4) above, or equivalently Theorem 3.2
(15) with  $m = 5$, occurs (actually realizable at least
for  $(g_2, ..., g_J) = (1, 1, 1)$);
see 4.6 for an alternative direct proof.
\roster

\medskip \noindent
We need the following result first.

\medskip \noindent
\plainproclaim Lemma 5.6. Let  $X$  be a Coble surface of rational type
with respect to a  $(-1)$-curve  $E$  and with  $\pi : X \rightarrow Y$  the blow-down
of
$E$.  Then we have:

\roster
\item
Suppose that  $X$  fits Case (14) (resp. Case (15)) of Theorem 3.2.  Then  
$Y = \pi(X)$  is basic if and only if  $m = 3$ (resp. $m = 4, 5$).

\item
If  $X$  fits Case (16), or Case (14) with  $k > 6$ ($=$ deg $(-2K_{\bbP^2})$), then
$Y$  is not basic.
\endroster

{\sl Proof.} Consider Theorem 3.2 (16).  The others are similar (see Remark 3.3 and the proof
of Lemma 4.2 for the ``if'' part of (1)).  In the following, we shall use  $M, G, H, H_1$
to denote their $\pi$-images on  $Y$ (cf. Lemma 2.3).

\medskip \noindent
\plainproclaim 5.6.1 Claim.  All negative curves ($\ne H_1$) 
on  $Y$  are contained in fibres of the  $\bbP^1$-fibration on  $Y$  induced 
from the one on  $Y_{\min}$.

\medskip
If the claim is false for some  $C$  on  $Y$, then  $C \cdot M \ge 1$
for  $\bar M^{\perp}_{\bbQ} = \bbQ \bar H_1$; using  $C$  to intersect
the equality  $-2K_Y = M + G + H$, we see that  $C$  is a  $(-1)$-curve
with  $(C \cdot M, C \cdot G + H) = (1, 1), (2, 0)$.
Expressing  ${\bar C} = \tau(C) \sim a \bar H_1 + b f$  with
a general fibre  $f$  on  $Y_{\min} = \bold F_m$,
we have  $b \ge am$  due to the irreducibility of  $\bar C$,
and get  $2 \ge C \cdot M = \bar C \cdot \bar M = b \ge am \ge m$,
a contradiction to the fact that  $m \ge 3$  in Theorem 3.2 (16).

\medskip
Let  $Y \rightarrow Y_1$  be the blow-down of all  $(-1)$-curves
in fibres disjoint from  $H_1, M$.  Then for each singular fibre  $F_i$  (of length
$n_i$) on  $Y_1$, the dual graph of  $H_1 + F_i + M$  on  $Y_1$  is as in
Lemma 1.8 with  $s_1 = H_1, s_2 = M$.
By the claim above, the basicness of  $Y$  would imply the existence of
a blow-down  $Y_1 \rightarrow \bold F_1$  of  $(-1)$-curves 
in fibres such that  $H_1$  becomes the unique  $(-1)$-curve on  $\bold F_1$;
hence if  $-b_1$  is the self-intersection of  $H_1$  on  $Y_1$, then
$\sum_i n_i \ge b_1 - 1$.  On the other hand, the intersection
of (the images of) $M$  and  $H_1$  on  $X, Y, Y_1, Y_{\min}$  are the same
by the construction of  $Y_1$  and by noting that  $M$  is the  $\tau$-pull back
of  $\bar M$  on  $Y_{\min}$ (Remark 3.3).  
So  $b_1 = b_1 + 2(M \cdot H_1) = M^2 + \sum_i n_i \ge m + b_1 - 1$
and  $m \le 1$ (cf. the proof of Lemma 1.8 and blow down  $Y_1$  to  
$\bold F_{b_1}$  to see the second equality).  
This contradicts the fact that  $m \ge 3$
in Theorem 3.2 (16).  So  $Y$  is not basic.

\medskip
{\bf 5.7}  Now we prove Theorem 5.4.  In view of Lemma 5.6, we only need
to consider Cases (12) and (13) of Theorem 3.2.
The former one will imply Theorem 5.4 (6) by the argument
in Lemma 5.6; indeed, all components of  $\pi(G + H)$  are
disjoint from  $C$  and all, except  $\pi(G_1), \pi(H_1)$,
contracted to points by the map  $Y \rightarrow \bold F_1$, while
$\pi(G_1), \pi(H_1)$ (resp. $\pi(M)$) are mapped to section(s) of 
self-intersection 1 (resp. $5$) on  $\bold F_1$ (cf. the proof of Lemma 1.8).

Consider Theorem 3.2 (13).  If all negative curves ($\ne \pi(G_1)$) on
$Y$  are contained in fibres, then the basicness of  $Y$  implies
that Theorem 5.4 (2) occurs.  Otherwise, the proof of Lemma 5.6
shows the existence of a section  $C$ ($\ne \pi(G_1)$)  on  $Y$  such that  $C^2 = -1$
and  $(k; \pi^{-1}(C) \cdot M_1, \pi^{-1}(C) \cdot G + H) =
(2; 1, 0), (1; 1, 1), (1; 2, 0)$.  In particular, $\pi^{-1}(C) \cdot G_1 = 0$
for  $G = 4G_1$  now.  The first two clearly imply Theorem 5.4 (7).

Now assume that  $(k; \pi^{-1}(C) \cdot M_1, \pi^{-1}(C) \cdot G + H) = (1; 2, 0)$.
We shall show that this will imply Theorem 5.4 (2).
Let  $Y \rightarrow Y_1$  be the blow-down of all  $(-1)$-curves in fibres 
disjoint from the double section  $C$.
Then for each singular fibre  $F_i$  on  $Y_1$,
either  $C + F_i$  is a simple loop so that  $F_i$  has
the same dual graph as its namesake in Lemma 1.8, 
or  $F_i = 2(E_i + H_i^{(1)} +
\cdots + H_i^{(n_i-2)}) + H_i^{(n_i-1)} + H_i^{(n_i)}$  where  $\sum H_i^{(j)}$
has type  $D_{n_i}$ Dynkin diagram ($n_i = 2, 3$  are possible),
where  $E_i$  is a  $(-1)$-curve meeting  $C$  and  $H_1$
(and also  $H_2$  when  $n_i = 2$).  Now utilizing the equality
$-2K_Y = \pi(kM_1 + G + H)$  and intersecting it with (inverses of) curves
in the fibre  $F_i$, we see that the loop case of  $C+F_i$  is impossible
and we have
$$-2K_{Y_1} = M_1 + 4G_1 + \sum_i [\sum_{j=1}^{n_i-2} 2j H_i^{(j)}
+ (n_i-2) H_i^{(n_i-1)} + n_i H_i^{(n_i)}], \eqno (5.1)$$
where we assume that the section  $G_1$  on  $Y_1$
meets the fibre  $F_i$  at  $H_i^{(n_i)}$.  
Intersecting  $G_1$  with  (5.1), one gets  $2G_1^2 = 3 - \sum_i n_i$.  
On the other hand,
the disjointness of the section  $G_1$  with the  $(-1)$-double section  
$C$  on  $Y_1$  implies that  $-4G_1^2 = C^2 + \sum_i n_i$ 
(cf. the proof of Lemma 1.8).
>From these two equalities, one deduces that  $G_1^2 = -1$  on  $Y_1$.  
Hence Case (2) occurs.  This proves Theorem 5.4.

\head 6. Rational curves with negative self-intersection \endhead
In this section we shall study $(-n)$-curves on a Coble surface. 
The goal is to see whether this set is finite, or
finite modulo automorphisms of the surface.  We start with a definition:

\medskip\noindent
{\bf 6.1}
Let $X$ be a Coble surface. We say that $X$ is of {\it $K3$-type} 
if $|-2K_X|$ contains a reduced divisor.  The reason
for this definition is explained by the folowing:

\medskip\noindent
\plainproclaim 6.2 Lemma. Let $X$ be a Coble surface. 
Then the following properties are equivalent:
\roster
\item  $|-2K_X|$
contains a reduced divisor.
\item  
There exists a double cover $\tilde X\to X$, where $\tilde X$ is 
a K3-surface with at most ordinary double points as singularities. 
\endroster

{\sl Proof.} (1) $\Rightarrow$ (2)
Let  $B \sim -2K_X$  be a a reduced effective anti-bicanonical
divisor.  Then  $B$  is of simple normal crossing (Lemma 1.4).
Let  $\tilde X$  be the double cover of $X$ corresponding to the square root of 
$B$  defined by the line bundle  $\calO_X(-K_X)$.  
By the formula for the canonical sheaf of a double cover
we get  $\omega_{\tilde X} = \calO_{\tilde X}$. 
Since  $B$  has at worst ordinary double points,
${\tilde X}$  is a K3 surface with at worst
ordinary double points.

(2)$\Rightarrow$ (1) This follows from the formula for the 
canonical class of a double cover.

\medskip\noindent
\plainproclaim 6.3 Theorem.  A Coble surface of rational type
with respect to some  $(-1)$-curve  $E$  will never be of K3-type.

{\sl Proof.} Suppose the contrary that  $X$  is a Coble surface
of rational type with respect to a  $(-1)$-curve  $E$, which is also of K3-type.  
So if  $\pi : X \rightarrow Y$  is the blow-down of  $E$,
then we have  $|-2K_X + 2E| = |M| + G + H = \pi^*(|-2K_Y|)$  with  $p_a(M) = 0$
and  $p_a(-2K_X+2E) = 1$.
By the condition and Lemma 2.3  to the extent that  $E \cap (G+H) = \emptyset$,
we see that  $G+H$  is reduced; in particular, the  $\tau$-image
$\Gamma = {\overline M} + {\overline G} + {\overline H}$
on  $Y_{\min}$  is also reduced.
By Remark 3.3 and calculating the image of  $M+G+H$  on  $Y_0$,
we see that only Cases (5), (6), (10), (11) are possible.

Assume Case (5) or (6) occurs and  $\Gamma = \bar M_1 + \sum_i \bar G_i$  
is of simple normal crossing;
the general case and Cases (10) and (11) are similar.
Noting that  $K_Y^2 \le 0$  and applying Lemma 1.9 repeatedly, 
we see that  $Y \rightarrow \bbP^2$  is the blow-up of the 9 
intersection points in  $\sum_i \bar G_i$;
we can not touch points on  $\bar M_1$
(see Remark 3.3).  Thus  $Y$  and  $X$  are equal to their namesakes 
in Example 2.12.  Hence  $X$  is of elliptic type with respect to  $E$,
a contradiction.  For general situation of Case (6) say, we need to apply
Lemma 2.2 (4) (the uniqueness of a loop, if exists, and the inequality there); 
in particular, all triple points as well as
all double points (with possibly one exception) of  $\sum \bar G_i$  must be
blown up; also note that there is no quadruple point of  $\sum \bar G_i$  due
to the reducedness of  $G + H$.  This proves Theorem 6.3.

\medskip
There is a strong relation between Coble surfaces of K3-type
and minimal resolutions of rational log Enriques surfaces of index 2.
A rational {\it log Enriques surface}  ${\bar X}$  of index 2  is 
a normal rational surface with at worst quotient singularities
such that  ${\Cal O}(-2K_{\bar X}) \cong {\Cal O}_{\bar X}$ 
(cf. [Zh1]).

\medskip\noindent
\plainproclaim 6.4 Proposition.

\roster
\item 
The minimal resolution  $X$  of a rational log Enriques surface
${\bar X}$  of index 2 is a Coble surface such that 
$h^0(-2K_X) = 1$  and the only member  $D$  in  $|-2K_X|$
is a reduced divisor whose connected component is
either a single  $(-4)$-curve or a linear chain with the
following dual graph:
$$(-3) - (-2) - \cdots - (-2) - (-3).$$
The converse is also true.
\item
A terminal Coble surface has exactly one
anti-bicanonical divisor  $D$, and  $D$
is reduced and a disjoint union of  $(-4)$-curves. 
The converse is also true.
\item
The minimal resolution  $X$  of a rational normal surface
${\bar X}$  with at worst type  $\frac{1}{4}(1,1)$
singularities is a terminal Coble surface.
The converse is also true.
\item
Let  $X$  be a Coble surface with a reduced divisor  $D \in |-2K_X|$.
Then there is an embedded resolution  $(X', D')$  of
$(X, D)$  with  $D'$  the proper inverse transform of  $D$,
such that  $X'$  is a terminal Coble surface with
$D'$  as the only member in  $|-2K_{X'}|$.
\endroster

{\sl Proof.} The first part of (1) is proved in [Zh1].
For the converse, if  $X \rightarrow \bar X$  is the
contraction of  $D$  then one sees easily that  $\bar X$
is a rational log Enriques surface of index 2.

We prove (2).  If  $X$  is terminal Coble, then an arbitrary
member  $D$  of  $|-2K_X|$  is reduced and smooth (Lemma 1.9)
and hence a disjoint union of  $(-n_i)$-curves  $D_i$ (Lemma 1.4).
Now  $D_i^2 = D \cdot D_i = D_i \cdot (-2K_X)$
implies that  $D_i$  is a  $(-4)$-curve; in particular,
$h^0(-2K_X) = h^0(D) = 1$.  This proves (2) (cf. Lemma 1.9).  

For the first part of (3), by the proof of (1),
$|-2K_X|$  has exactly one member  $D$
which is reduced and a disjoint union of  $(-4)$-curves.
So  $X$  is a terminal Coble surface (cf. Lemma 1.9).
For the converse of (3),
we let  $X \rightarrow \bar X$  be the contraction of the unique
divisor  $D$  in  $|-2K_X|$.
Then  $\bar X$  satisfies the required condition.

Next we prove (4).  By Lemma 1.4, $D$  has only nodes 
as singularities.  Let  $X' \rightarrow X$  be the blow-up
of all nodes in  $D$.  Then we have  $-2K_{X'} \sim D'$.
This implies, as in (2), that  $D'$  is a disjoint
union of  $(-4)$-curves.  Hence  $X'$  is terminal Coble.
This completes the proof of Proposition 6.4.

\medskip
Let  $f:X' \to X$ be  a birational morphism of Coble surfaces. 
If  $X'$  is of  $K3$-type then so is  $X$;
indeed, if  $D' \in |-2K_{X'}|$  is reduced then so 
is  $f_*(D') \in |-2K_X|$.
In view of the above observation and Lemma 1.10,
among Coble surfaces of K3-type, minimal ones are the most
interesting.  Such  $X$  is of elliptic type with respect to any
$(-1)$-curve  $E$ (Theorem 6.3).  Suppose that
$M^2 = 0$  in notation of Lemma 2.2.  
Then  $X$  is given in either Theorem 2.5 or Theorem 2.8
with  $X = X'$  and  $E = E'$.

\medskip \noindent
\plainproclaim 6.5 Theorem.  Suppose  $X$  is a Coble surface with  $M^2 = 0$.  
If  $X$  is of Halphen type obtained from a minimal 
Halphen surface  $Y_m$  of index 2 by one blow-up of 
a singular point on its non-multiple fibre  $F$, 
then it is of K3-type if and only if  
$F$  is  of type  $I_n$, II, III or IV. 
If  $X$ is of Jacobian type obtained as in Theorem 2.8
from a minimal Jacobian rational elliptic surface  $Y_{\min}$
by blowing up a singular point from one fibre  $F$  and
singular points (at least one) and their infinitely near points
on another fibre  $F_1$,
then it is of K3-type if and only if each of  $F$  and  $F_1$  is
of type $I_n$, II, III, or IV.

{\sl Proof.} This follows immediately from the Kodaira classification of singular fibres.

\medskip
There is an analogue of the K3-cover for Coble surfaces
of elliptic type which are not of K3-type:

\medskip\noindent
\plainproclaim 6.6 Theorem. Suppose  $X$  is a Coble surface of elliptic
type with  $M^2 = 0$  in notation of 2.1, which is not of K3-type.  
Then  $X$  admits a double cover  $\tilde X$  which is
a non-minimal rational Jacobian elliptic surface.

{\sl Proof.} We do only the case when  $X$  is of Halphen type; 
the Jacobian case can be considered similarly.
Then  $X$  is obtained from a minimal Halphen elliptic surface  $V$  of index 2
by blowing up a singular point of its non-multiple fibre  $F$ 
of type  $\ne I_n,II,III, IV$.

We check the assertion by considering different types of the
fibre.  Let us do for example, the case  $F$  is of type $I_{b}^*$ and
leave the other cases to the reader. 
Write  $F = R_1+R_2+R_3+R_4+2(R_5+\ldots R_{b+5})$, 
where  $R_1,R_2$  intersect  $R_5$
and  $R_3, R_4$  intersect  $R_{b+5}$. Then 
$$R_1+R_2+R_3+R_4 \sim -2K_X-2(R_5+\ldots R_{b+5}),$$
hence there exists a double cover $\pi:\tilde V \to V$ ramified over
$R_1+R_2+R_3+R_4$.  We have 
$$K_{\tilde V} \sim - \pi^{-1}(R_5+\ldots R_{b+5}).$$
If  $b = 0$, $C = \pi^{-1}(R_5)$  is an elliptic curve with  $C^2 = -4$. 
If  $b\ne 0$, $C = \pi^{-1}(R_5+\ldots R_{b+5})$  is a reducible curve of
arithmetic genus 1. The pre-image of a general fibre of the elliptic
fibration on  $\tilde V$  splits into a disjoint union of two elliptic curves. 
After a base change $\bbP^1\to \bbP^1$ of degree 2 ramified at two points, we obtain an
elliptic fibration on 
$\tilde V$  with one of its fibre equal to  
$(\pi^*(F))_{\text{\rm red}} = C+ (\tilde R_1+\tilde R_2+\tilde R_3+\tilde R_4)$,
where  $\pi^*(R_i) = 2\tilde R_i, 1 \le i \le 4$.  
Note that the  $\tilde R_i$  are  $(-1)$-curves on  $\tilde V$.  
Blowing these four curves down, we obtain an elliptic surface  $\hat V$  with the
image  $\hat C$  of  $C$, which is reduced and linearly equivalent to $-K_{\hat V}$. 
One can verify that  $\hat V$  is a Jacobian Halphen surface.

Now, if  $X$  is obtained from  $V$  by blowing up a
point  $p$  on $R_5+\ldots + R_{b+5}$, it admits a double cover  $\tilde X$
which is obtained from  $\tilde V$  by blowing up two (or one if $p$
also lies on some  $R_i$  with  $i\le 4$) points on
$C$. So $\tilde X$ is obtained from the minimal elliptic surface 
$\hat V$  by blowing up points on one fibre $\hat C$. 

\medskip
Now we consider the finiteness problem 
of the number of negative curves on a Coble surface modulo automorphisms.

\medskip\noindent
\plainproclaim 6.7 Theorem. Assume $k= \bbC$. Let $X$ be a Coble surface
of elliptic type. Suppose that $X$ is a terminal Coble surface of K3-type. Also
assume that  $X$ is general in the sense that any divisor class on the K3-cover
is invariant with respect to the double cover involution. Then the group
$\roman{Aut}(X)$ has finitely many orbits in the set of negative rational
curves on $X$.

{\sl Proof.} This follows from two well-known results about K3-surfaces.
The first one says that the group of automorphisms of any K3-surface has
only finitely many orbits in the set of smooth rational curves (see
[Na, St]). The second one says that any automorphisms of the K3-cover
of $X$ commutes with the involution (see [Ni]).

\medskip
We do not know whether the same result is true for non-terminal Coble surfaces of
K3-type. However we shall show now that it cannot be extended to Coble surfaces not
of K3-type.

\plainproclaim  6.8 Lemma.
Let  $\pi_A : S(A) \to S$  be the blow-up of a set  $A$  of $n$ points
on a nonsingular projective surface $S$ with zero irregularity.  
Let  $G(A)$  be the subgroup of $\roman{Aut}(S(A))$  consisting of 
automorphims which are identical on the proper inverse transform $C'$ 
of a nonsingular irreducible curve $C$ of positive genus on
$S$  which contains $A$. Then the set of subsets $A$ of $C$ such that
$G(A)$ is not the lift of a subgroup
$G(A)'$ of $\roman{Aut}(S)$ is countable.

{\sl Proof.} We use induction on $n$. Assume $n = 1$. Let  $E_A$  be
the exceptional curve of  $\pi_A$.  An element
$g\in G(A)$ is a lift of an automorphism of $S$ if and only if $g$
stabilizes  $E_A$. Suppose $g(E_A)\ne E_A$.  The image 
$R_A$  in  $V$  of  $g(E_A)$  intersects  $C$  at one point
$a$  with multiplicity  $m+1$, where  $m = E_A\cdot g(E_A)$. 
The restriction of the linear system  $|R_A|$  to  $C$  is of degree
$m+1$, so that there are only finitely many points $c$ on  $C$  
which can be realized as a divisor $(m+1)c$ from this linear system. 
Here we use that the Jacobian of a curve of positive genus has only 
finitely many points of given finite order.  Since the set of
divisor classes on a surface with zero irregularity is countable, 
only a countable set of points $a\in C$ may have the property 
$g(E_A)\ne E_A$. This proves the assertion for $n = 1$.

If $n > 1$, we write $A = A'\cup \{a\}$, where $a\not\in A'$. The map
$\pi_A$ is equal to the composition of the maps 
$\pi_a:S(A)\to S(A')$ and $\pi_{A'}:S(A')\to S$. It is clear that $C'$ is
equal to the proper inverse transform of a curve $C''$ on $S(A')$ 
which is also the proper inverse transform of $C$. By the case $n = 1$, 
we know that the set of points $a$ for which elements of $G(A)$ do not
descend to automorphims of $S(A')$ is countable. By induction, the set of
subsets $A'$ for which elements of $G(A)$ do not descend further to $S$ 
is countable. So, the set of all
possible $A$ for which elements of $G(A)$ do not descend to $S$ is countable. 

\medskip
\plainproclaim 6.9 Lemma. Let  $\Sigma$ be a set of 9 points in
$\bbP^2$  and let  $X$  be the blow-up of  $\Sigma$.  Denote by  $\calE$  the set of 
all  $(-1)$-curves on  $X$.  Assume that  $\calE$  is infinite. 
Then, for any $E\in \calE$, the image  $S$  of the map 
$\calE \to \bbZ$  given by  $E' \to E' \cdot E$, is an infinite set.

{\sl Proof.} We have  $(E'-E)^2 = -2-2E' \cdot E$, so it suffices to show that
the set $S'$ of possible integers $m$ of the form $m = (E'-E)^2$ is infinite.
Since  $(E'-E)\cdot K_X = 0$, the divisor class of $E'-E$  belongs to the orthogonal
complement  $(\bbZ K_X)_{\Pic(X)}^\perp$. Since $K_X^2 = 0$, the
lattice $L = (\bbZ K_X)_{\Pic(X)}^\perp/\bbZ K_X$ is negative definite. This implies
that the set of vectors in $L$ of fixed norm is a finite set. In particular, if $S'$
is finite, the set of cosets in  $L$  of the classes $E'-E$ is finite.
On the other hand, $(E'-E)-(E''-E) = E'-E''\in \bbZ K_X$
would imply that  $0 = (E'-E'')^2 = -2-2E'\cdot E''$  and hence  $E' = E''$.
So  $E'-E$  and  $E''-E$ cannot belong to the same coset modulo  $\bbZ K_X$ 
unless $E' = E''$.  This shows that the set of divisor classes of $E'-E$ must be
finite, contradicting the assumption that $\calE$ is infinite. 
This contradiction proves the lemma.

\medskip\noindent\medskip\noindent
{\bf 6.10 Example}. Here we give an example of a minimal Coble surface of
Halphen type such that its automorphism group has infinitely many
orbits on the  set of $(-1)$-curves. We have to assume that the ground field $k$ is
uncountable.

Let $V$ be an Halphen surface
of index 2 with a reducible fibre of type $I_0^*$.  One can explicitly
construct it as follows. Take five lines $L_i$ ($1 \le i \le 5$) in $\bbP^2$ in
general linear position and consider a pencil of elliptic curves spanned by the
curve $C_6 = L_1+L_2+L_3+L_4+2L_5$ and the curve $2C_3$, where $C_3$ is
the cubic which passes through 6 intersection points $p_{ij} = L_i\cap L_j;
i,j = 1,\ldots,4$. We assume that the cubic $C_3$ intersects $L_5$ at
three distinct points $q_1,q_2,q_3$. Resolving the base points of the
pencil we arrive at a Halphen surface $V$ of index 2. The image to
$\bbP^2$ of its fibre of type $I_0^*$  is the sextic $C_6$.
If we choose a  point $a\in L_5$ and blow up the corresponding
point on $V$ we obtain a minimal Coble surface $X$ of Halphen type.

Let $V'$ be the Halphen surface obtained in the same way as $V$
but replacing the cubic curve by a new cubic which passes through the
points $p_{ij}$ and the points $q_1,q_2, a$. We have a natural map $f : X \to V'$
which is the blow-up of the pre-image $q_3'$ of $q_3$ on $V'$. When
$a$ is chosen general enough, the elliptic fibration $\pi:V'\to \bbP^1$ has only
one reducible fibre (of type $I_0^*$). Its Jacobian fibration (a relative
minimal model of the Jacobian of the generic fibre of $\pi$) has only one
reducible fibre of type $I_0^*$ and hence its Mordell-Weil group $\roman{MW}$ is
infinite (of rank 4). Since $\roman{MW}$ acts freely
by translations on the set of bi-sections of $\pi$, we see that $V'$ has
infinitely many $(-1)$-curves (which are rational bi-sections of $\pi$). 
Let  $E_a$  be the exceptional curve on  $V'$  blown-up from the point $a$. 
Its pre-image, also denoted by  $E_a$, under the map  
$f: X \to V'$, is the exceptional curve of 
the map $X\to V$. By Lemma 6.9, $V'$  has $(-1)$-curves  $E_i$  with unbounded set of
integers  $m_i = E_i \cdot E_a$.  The pull-backs on  $X$, also denoted by  $E_i$,
of the curves  $E_i$  on $V'$, form an infinite set of $(-1)$-curves (if $E_i$ does not pass
through $q_3'$) or $(-2)$-curves (if $E_i$ passes
through $q_3'$) with unbounded intersection numbers with a general fibre of
the elliptic fibration  on  $X$ (the pull-back of the elliptic fibration on  $V$). 
Since the set of $(-1)$-curves on $V'$ is countable we can always choose $a$ and
$q_3$ such that $(-1)$-curves on $V'$ do not pass through $q_3'$. So we can assume that
all  $E_i$'s are $(-1)$-curves.  Thus we have found infinitely many 
$(-1)$-curves $E_i$ on $X$ with unbounded intersection numbers with 
a general fibre of the elliptic fibration  on  $X$.

Note that an automorphism  $g$  of  $X$  leaves invariant the 
isolated linear system  $|-2K_X| = \{ R_1 + \cdots +  R_4 + 2 R_5\}$,
where $ R_i$ denotes the proper inverse transform of $L_i$ in $X$. In particular,
$ R_5$  is  $g$-stable. Let  $G$  be the kernel of the natural action of  
$\roman{Aut}(X)$  on the 4-point set  $\{ R_1, \cdots,  R_4\}$.  Then
$\roman{Aut}(X)/G$  is isomorphic to a subgroup of the
symmetric group  $S_4$  in 4 letters.
Now each  $g \in G$  fixes all 4 points  $ R_i \cap R_5$
of the rational curve  $R_5$  and hence  $g$  acts identically on  $R_5$
(there is no non-trivial automorphism of $\bbP^1$  which
fixes more than two distinct points).
Take the double cover $S \to V$ branched along the union of the curves $R_1+\ldots+R_4$ (see
Theorem 6.6). The pre-image of $R_5$ is an elliptic curve $C$ on $S$. Let $A = \{a',a''\}$ be
the pre-image of $a\in R_5$ on $C$. Consider the group $G(A)'$ of automorphims of the blow-up
$S(A)$ of $A$ which are lifts of automorphisms $g\in G$. Recall that,
since all elements of $G$ leave the square root invariant of 
the divisor class of the branch divisor, 
for every $g\in G$ there is an element 
$\tilde g \in \roman{Aut}(S(A))$ which commutes with the involution $\sigma$ of 
the double cover, and descends to an automorphism of $X$. Two lifts of the 
same $g$ differ by $\sigma$.   
All elements of $G(A)'$ restrict to $C$ as automorphims of order $\le 2$. 
Let $G(A)$ be the subgroup of index
$2$ of $G(A)'$ consisting of elements of $G(A)'$ which act identically on $C'$. 
We shall identify this group with the group $G$. By Lemma
6.8, we can choose $a$ such that all elements of $G(A)$ are lifts of automorphisms of
$S$. Thus all elements of $G$ are lifts of automorphims of $V$ to $X$. In particular,
$g$  stabilizes  $E_a$  so that the full fibre  $R_1 + \cdots + R_4 + 2(R_5 + E_a)$
on  $X$  is  $g$-stable.
Clearly each  $g \in G$ preserves the degrees of the multi-sections
$E_i$.  Hence the number of orbits of  $G$ (and also of $\roman{Aut}(X)$,
due to the finiteness of the index of  $G$  in it) 
on the set of $(-1)$-curves on  $X$  is infinite.

\medskip\noindent 
{\bf 6.11}  Next we would like to study the set $\calE$ of $(-n)$-curves 
($n \ge 1$) on a Coble surface  $X$  of rational type. 
We still do not have a
complete picture of $\calE$; however we guess that this set is always
finite.  In fact, by a theorem of Nagata ([Na], Theorem 5) this is always
true if  $X$ is not basic.  Another special case where it is true is when
$\kappa^{-1}(X) = 2$. 
Here  $\kappa^{-1}(X)$  denotes the {\it anti-Kodaira} dimension. 
This is the Iitaka-Kodaira dimension of the divisor $-K_X$.  
Obviously  $\kappa^{-1}(X)\ge 0$ for Coble surfaces.  It is also clear that
$\kappa^{-1}(X) \le 1$ for Coble surfaces  $X$  of elliptic type with 
$M^2 = 0$  in notation of Lemma 2.2.

We shall use the following result from [Sa]:

\medskip\noindent
\plainproclaim  Lemma. Let $X$ be a surface with $\kappa^{-1}(X) = 2$. 
Then  $X$ has only finitely many curves with negative self-intersection. 

\medskip\noindent 
{\bf 6.12} Let $X$ be a Coble surface of rational or elliptic type with
respect to a curve $E$ and let  $\pi : X \rightarrow Y$ 
be the blow down of  $E$.  By the definition,
$\pi^*(|-2K_Y|) = |M|+P$, where $M^2 \ge 0$.  We have
$\kappa^{-1}(Y) = 2$  if and only if either  $M^2 > 0$, or  $M^2 = 0$  and
$p_a(M) = 0$  (noting that then  $P \cdot M = 4k > 0$, see Lemma 2.2);
if this is the case, $Y$ contains only finitely many negative rational curves. 
Unfortunately, this does not automatically imply the finiteness of $\calE$ for $X$
except in a few special cases which we list now.  

\medskip\noindent 
\plainproclaim 6.13 Lemma. Let $f:X'\to X$ be the blow-up of a point
$p\in X$. Then  $\kappa^{-1}(X') \le \kappa^{-1}(X)$. 
The equality takes place if $p$ is a point
of multiplicity $\ge {n+1}$ of an effective divisor from $|-nK_X|$. 

{\sl Proof.}  The first assertion is obvious. The second assertion follows from the
fact that the anti-Kodaira dimension of a divisor $D$ depends only on
$D_{\roman{red}}$. 

\medskip\noindent 
\plainproclaim 6.14 Proposition. In notation of Lemma 2.2, 
\roster
\item
Assume that  $M' \ge 3E$  for some  $M'$  in  $|M|$
and  $(p_a(M), M^2) \ne (1, 0)$.  Then  $\kappa^{-1}(X) = 2$.
\item
If  $M = kM_1$  with  $k \ge 3$.  Then  $\kappa^{-1}(X) = 2$.
\item
One has  $\kappa^{-1}(X) = 2$  if either  $p_a(M) = 1$
and  $M^2 = 6$, or  $p_a(M) = 0$  and  $M^2 \ge 4$.  
\endroster

{\sl Proof.} For (1), we consider only the case where  $M^2 = 0$  and  $p_a(M) = 0$.
Then  $P \cdot M_1 = 4$ (Lemma 2.2).  Thus  
$\kappa^{-1}(X) = \kappa(X, M_1+P) = 2$, where the first equality follows
from the observation that  $-2K_X \sim (M'-2E) + P$  and the latter
has the same support as  $M' + P$ ($\sim kM_1 + P$).  (2) is a consequence of (1).

For (3), we consider only the case  $p_a(M) = 0$.
By Lemma 1.7, $h^0(M) = M^2 + 2$.  Hence if  $M^2 \ge 5$, then
there is a member  $M'$  in  $|M|$  with  $M' \ge 3E$, or
equivalently  $\pi(M')$  has multiplicity $\ge 3$  at the point
$\pi(E)$ (cf. Lemma 2.3).  So (3) is true in this case.
Suppose that  $m = M^2 = 4$.  Then Case (9), (14), (15) or (16) in Theorem 3.2
occurs.  We treat Case (15) because the others are similar.
Now  $M$  is the $\tau$-pull back of  ${\bar M} \sim \bar G_1 + (m-2)f$,
where  $f$  is a fibre and  ${\bar G}_1$  a section of self-intersection  $-(m-4) = 0$
(cf. Remark 3.3).  Let  $M'$  be the sum of the  $\tau$-pull backs of
$\bar G_1' \in |\bar G_1|$  and  $2f' \in |2f|$
through the point  $\tau(E)$  and (1) applies.
One can actually shows that (3) is still true even when  $p_a(M) = 0$
and  $M^2 = 3$, unless Theorem 3.2 (14) with  $m = 3$  occurs.

\vfill\eject

\vglue 1in
\noindent
Igor V. Dolgachev: Department of Mathematics, University of Michigan
\smallskip\noindent
Ann Arbor, MI 48109, USA
\smallskip\noindent
e-mail: idolga$\@$math.lsa.umich.edu

\bigskip\noindent
De-Qi Zhang: Department of Mathematics, National University of Singapore
\smallskip\noindent
Kent Ridge, Singapore
\smallskip\noindent
e-mail: matzdq$\@$math.nus.edu.sg

\vfill\eject
\bye